\documentclass{amsart}
\usepackage{amsmath}
  \usepackage{paralist}
  \usepackage{graphics} 
  \usepackage{epsfig} 
 \usepackage[colorlinks=true]{hyperref}
\hypersetup{urlcolor=blue, citecolor=red}

\newtheorem{theorem}{Theorem}[section]

\newtheorem{lemma}[theorem]{Lemma}

\newcommand{\sgn}{\mbox{sgn}}
\newcommand{\vct}[1]{{\mathbf #1}}       
\newcommand{\mtx}[1]{{\mathrm #1}}       

\newcommand{\SSS}{{\mathcal{S}}}

\newcommand{\bound}{\partial\Omega}

\newcommand{\abs}[1]{\left\vert {#1} \right\vert}
\newcommand{\norm}[1]{\left\Vert {#1} \right\Vert}

\newcommand{\pd}[2]{\frac{\partial {#1}}{\partial {#2}}}

\newcommand{\im}{\mbox{Im\,}}
\newcommand{\re}{\mbox{Re\,}}
\newcommand{\dbar}{\bar{\partial}}

\newcommand{\CC}{{\mathcal{C}}}

\newcommand{\OO}{{\mathcal{O}}}

\newcommand{\A}{{\mathcal{A}}}

\newcommand{\T}{{\mathbf{t}}}
\newcommand{\Om}{\Omega}
\newcommand{\DOm}{\partial\Omega}

\newcommand{\R}{{\mathbb R}}

\newcommand{\C}{{\mathbb C}}

\newcommand{\matr}[1]{\left[\!\begin{array}{l} #1 \end{array}\!\right]}

\newcommand{\kk}[1]{{\overline{#1}}}

\newcommand{\TT}{{\mathcal{T}}}

\newcommand{\MHz}{\textrm{ MHz}}
\newcommand{\supp}{\textrm{supp}}
\newcommand{\dee}{{\rm d}}
\newcommand{\rmi}{\mathrm{i}}

\newcommand{\imagepath}{}
\newcommand{\refpath}{}

\title[The D-bar Method for Diffuse Optical Tomography]{The D-bar Method for Diffuse Optical Tomography: a computational study}

\author{J.~P.~Tamminen, T.~Tarvainen and S.~Siltanen}


 \email{janne.tamminen@ttu.ee}
 \email{tanja.tarvainen@uef.fi}
 \email{samuli.siltanen@helsinki.fi}


\begin{document}

\begin{abstract}
The D-bar method at negative energy is numerically implemented. Using the method we are able to numerically reconstruct potentials and investigate exceptional points at negative energy. Subsequently, applying the method to Diffusive Optical Tomography, a new way of reconstructing the diffusion coefficient from the associated Complex Geometrics Optics solution is suggested and numerically validated.
\end{abstract}

\maketitle



\bigskip

\tableofcontents

\section{Introduction}
\noindent
Diffuse optical tomography (DOT) is an imaging modality in which images of the optical properties of the medium are estimated based on measurements of near-infrared light on the surface of the object. In practice, light is guided to the surface of the subject, and the transilluminated light measured on the surface, using fibre optics. DOT has potential applications in medical imaging, for example in breast cancer detection and classification, monitoring of infant brain tissue oxygenation level and functional brain activation studies, for reviews see e.g. \cite{Arridge1999,Gibson2005}.

The image reconstruction task of DOT is an example of an {\em ill-posed inverse problem.} This means that even quite different targets may produce almost the same data, leading to straightforward numerical inversion to be highly sensitivite to modelling errors and measurement noise. Therefore, any robust reconstruction method for DOT needs to be {\em regularized.} In this work we propose a non-iterative DOT reconstruction method where regularization is provided by a low-pass filter applied in nonlinear Fourier transform domain: we develop the numerical aspects of \textit{the D-bar method at negative energy}, in dimension two.
  
Absolute imaging in DOT uses a single set of measurements to reconstruct the spatial distribution of optical parameters  \cite{Arridge1999}. 
Generally, this image reconstruction problem is formulated as a numerical minimisation problem such as regularised least squares. The minimisation problem is solved using some iterative algorithm such as conjugate gradients or Gauss-Newton method.  The solution of this minimisation problem can be time consuming since, in order to find the minimum between measurements and model predictions, one needs to solve the forward model, the diffusion equation, in a discretisation covering the whole target domain at each iteration. This can be both memory and time consuming.

The D-bar method is a direct reconstruction method, meaning that the unknown is linked to the ideal measurement directly via certain equations, and the forward model needs not to be solved, except in simulations of data. The D-bar method brings the benefits of a direct reconstruction algorithm to DOT including the absence of local minima and the ability to parallelise the reconstruction computation. Furthermore, the approach can be used to reconstruct only the region of interest. It can also be used to compute an initial guess for other methods. 

In subsection \ref{sec:Dbar} we briefly present the D-bar method, however for more information and theoretical background, we refer to the review \cite{Grinevich2000}, the fundamental papers \cite{Novikov1992}, \cite{Grinevich1988} and the references therein. The original method described in \cite{Grinevich2000,Novikov1992,Grinevich1988} has recently been generalized in \cite{Lakshtanov2015b} and \cite{Lakshtanov2015c}, to handle the presence of \textit{exceptional points}, which traditionally have prevented the use of the D-bar method. However, in our work we concentrate on the original method and leave the application of these new results to future works.



For the remainder of the paper, we refer to the original D-bar method as \textbf{the D-bar method.} The main goals and novelties of this paper are:
\begin{itemize}
	\item numerically implement and test the D-bar method at negative energy, see section \ref{sec:Dbarnumerics},
	\item numerically find exceptional points which prohibit the use of D-bar, see section \ref{sec:excep},
	\item introduce and numerically validate a new reconstruction step for a \textit{conductivity at negative energy}, see section \ref{sec:conjecture}, and
	\item numerically test the D-bar method in a DOT setting, see section \ref{sec:DOT}.
\end{itemize}

The D-bar method is a convenient approach to nonlinear inverse problems because it directly reconstructs the desired PDE coefficient without iterations. However, the application to DOT has the theoretical limitation in that the nonlinear Fourier transform may be singular in some cases, preventing stable inversion. More precisely, the nonlinear Fourier transform uses \textit{complex geometrical optics} (CGO) solutions in place of exponential functions, and the existence and uniqueness of CGO solutions is not completely understood at present. Our numerical experiments shed light on this issue, showing two new things: (1) some practically relevant PDE coefficients have unique CGO solutions and allow reconstruction using the D-bar method; and (2) it is possible to construct coefficients that do exhibit singularities in their nonlinear Fourier transform. This is the first work that computationally explores such singularities at negative energy, see section \ref{sec:excep}.

Also other direct inversion methods for DOT have been studied. These are typically associated for particular experimental geometries such as  planar, cylindrical or spherical boundaries \cite{Arridge2009}. Often, a linearised inverse problem is considered. However, extensions to more complex geometries \cite{Ripoll2006} as well as non-linear approaches have been investigated \cite{Moskow2008,Moskow2009}.

The introduction continues with some mathematical background. Throughout the paper we identify the plane $\R^2$ as the complex plane by writing 
$$
z = [x_1\ x_2]^T\in\R^2,\quad z = x_1+\rmi x_2\in\C.
$$

\subsection{The Gel'fand-Calderon problem}
The origin of the D-bar method lies in the solution of the \textit{Gel'fand-Calderon problem} (GC problem) \cite{Gelfand1961,Calder'on1980}. Let $\Omega\subset\R^2$ be the unit disk. Let $q=q_0-E$ be the potential with the energy $E\in\C$ and $\supp(q_0)\subset\Omega$. Consider the Schr\"odinger equation
\begin{equation}\label{schrodinger-nonzero}
(-\Delta+q)u = 0 \quad \textrm{in }\Omega
\end{equation}
with the boundary condition $u=f$ on $\DOm$. This boundary-value-problem is well-posed if zero is not a Dirichlet-eigenvalue of $-\Delta+q$, then for any $f\in H^{1/2}(\DOm)$ it has a unique weak solution $u \in H^1(\Omega)$. For well-posed problems we define the \textit{Dirichlet-to-Neumann-map} (DN-map)
\begin{equation} \label{DNq-def}
\Lambda_q:\quad H^{1/2}(\DOm)\rightarrow H^{-1/2}(\DOm),\quad f\mapsto \pd{u}{\nu}\rvert_{\DOm},  
\end{equation}
where $\nu$ is the unit outer normal to the boundary. The inverse problem of our interest is then the following: given $\Lambda_q$ and the energy $E$, reconstruct the potential $q_0$.


The GC problem is a way to formulate different tomographic methods. Electrical Impedance Tomography (EIT), Acoustic Tomography (AT) and DOT are related to the GC problem by a transformation resulting to different energies. EIT is a zero-energy problem with $E=0$, AT is a positive energy problem with $E>0$ and DOT is a negative energy problem $E<0$.



\subsection{Diffuse Optical Tomography}
Consider the diffusion approximation problem
\begin{equation}\label{DOT}
\left\{
\begin{array}{rcll} 
-\nabla\cdot D\nabla \tilde{u} +(\mu_a+\rmi\frac{\omega}{c})\tilde{u} &=& 0    &\textrm{ in }\Om\\
                                  \tilde{u} + 2D\pd{\tilde{u}}{\nu}&=& g^-  & \textrm{ on }\DOm \\
                                               D\pd{\tilde{u}}{\nu}&=& -g^+ & \textrm{ on }\DOm,
\end{array}\right.
\end{equation}
where $\tilde{u}$ is the light fluence rate, $D$ is the diffusion coefficient, $c$ is the speed of light in the medium, $\omega$ is the angular modulation of the input, $\mu_a$ is the absorption coefficient, $g^-$ is the source on the boundary and $g^+$ is the scattered field measured on the boundary. See \cite{Arridge1999} for a review of DOT and the diffusion approximation. Define the DN-map by
\begin{equation*} \label{DN-DOT}
  \Lambda_{D,\mu_a}(g^-+2g^+) = -g^+.
\end{equation*}
The inverse problem of DOT is to reconstruct the diffusion coefficient $D$ and the absorption coefficient $\mu_a$ from the knowledge of the DN-map $\Lambda_{D,\mu_a}$.

DOT is related to the GC problem at negative energy in the following way. Assume
\begin{equation}\label{boundary-values}
D|_{\DOm} = d>0,\quad \mu_a|_{\DOm} = m>0.
\end{equation}
Then by writing $u = D^{1/2}\tilde{u}$ we get the Schr\"odinger equation \eqref{schrodinger-nonzero} of the GC problem with 
\begin{equation*}
q_0 = D^{-1/2}\Delta D^{1/2}+\frac{1}{D}(\mu_a+\frac{\rmi\omega}{c})+E,\quad E=-(\frac{m}{d}+\frac{\rmi\omega}{dc})\in\C
\end{equation*}
The DN-maps for the two problems will be related by
\begin{equation*}
\Lambda_q = \frac{1}{d}\Lambda_{D,\mu_a}.
\end{equation*}
To simplify the problem we can assume that $q$ is real-valued, then
\begin{equation*} \label{DOT-potential}
q_0 = D^{-1/2}\Delta D^{1/2}+\frac{\mu_a}{D}+E,\quad E=-\frac{m}{d}.
\end{equation*}
With the D-bar method we can reconstruct $q_0$ which includes both $D$ and $\mu_a$. Note that the first term of $q_0$ is a potential of \textit{conductivity type} which was originally defined by Nachman in the EIT problem \cite{Nachman1996}.




\subsection{The D-bar method}\label{sec:Dbar}
In short, the D-bar method is based on a non-linear Fourier transform where exponentially behaving CGO solutions of Faddeev \cite{Faddeev1966} are used. The boundary integral equation proved by R.~G.~Novikov \cite{Novikov1988} links the measurements to these CGO solutions. The D-bar equation discovered by Ablowitz, Nachman, Beals and Coifman \cite{Beals1981,Ablowitz1986} reveals the pseudoanalytic nature of the CGO solutions. The method works in the abscence of exceptional points for which the unique CGO solution does not exist. Recently in the works of Novikov and E.~Lakshtanov \cite{Lakshtanov2015b,Lakshtanov2015c} the D-bar method has been generalized to handle exceptional points, however we do not discuss these latest results here but instead concentrate on the original method.

We define Complex Geometrics Optics solutions as exponentially growing solutions $\psi(\,\cdot\,,\zeta)$ of the Schr{\"o}dinger equation \eqref{schrodinger-nonzero} in the whole plane $\R^2$, where $\zeta=[\zeta_1 ,\\ \zeta_2]^T\in\C^2$ is a spectral parameter with $\im(\zeta)\neq\vct{0}$. For a given potential $q_0$ and parameter $\zeta$ there might not exist a unique CGO solution in which case we call $\zeta$ an exceptional point.

Let us parametrize the subset of parameters $\zeta$ satisfying $\zeta\cdot\zeta=E$ by
\begin{equation}\label{lambda-definition}
\lambda = \frac{\zeta_1+\rmi\zeta_2}{\sqrt{E}},\quad \zeta = \matr{(\lambda+\frac{1}{\lambda})\frac{\sqrt{E}}{2} \\ (\frac{1}{\lambda}-\lambda)\frac{\rmi\sqrt{E}}{2}}. 
\end{equation}
We call the parameter $\lambda$ the spectral parameter as well. We write
$$
\mu(z,\zeta) = e^{-\rmi\zeta\cdot z}\psi(z,\zeta)
$$
and call it a CGO solution as well. Depending on whether we use $\zeta$- or $\lambda$-notation, in place of $\psi(z,\zeta)$ and $\mu(z,\zeta)$ we write $\psi(z,\lambda)$ and $\mu(z,\lambda)$ respectively, even if the energy is then omitted. The CGO solution $\mu(z,\lambda)$ satisfies
\begin{equation}\label{Llambda-discrete}
  (L_\lambda+q_0)\mu(z,\lambda) = 0, \qquad L_\lambda = -4\partial_z\dbar_z-2\rmi\sqrt{E}(\lambda\partial_z+\frac{1}{\lambda}\dbar_z).
\end{equation}
The Green's function of $L_\lambda$ is called \textit{Faddeev Green's function} and it is denoted as $g_\lambda(z)$. The CGO solution satisfies the Lippman-Schwinger (LS) equation
\begin{equation}\label{LS-lambda}
  \mu(z,\lambda)= 1-g_\lambda(\cdot)\ast(q_0(\cdot)\mu(\cdot,\lambda))(z).
\end{equation}

In order to describe how the D-bar method solves the inverse problem, some more machinery is needed. Define the differential operators
\begin{equation*}
\partial_w = \frac{1}{2}(\partial_{w_1}-\rmi\partial_{w_2}),\quad \dbar_w = \frac{1}{2}(\partial_{w_1}+\rmi\partial_{w_2}),
\end{equation*}
where $w = w_1+\rmi w_2$, and the exponential functions
\begin{eqnarray}
e_{-\lambda}(z) &=& \exp\left(-\frac{\rmi\sqrt{E}}{2}(1-\frac{1}{\lambda\kk{\lambda}})(-z\kk{\lambda}+\kk{z}\lambda)\right)\label{exp-mlambda},\\
e_{\lambda}(z)  &=& \exp\left(\frac{\rmi\sqrt{E}}{2}(1-\frac{1}{\lambda\kk{\lambda}})(-z\kk{\lambda}+\kk{z}\lambda)\right)\label{exp-lambda}.
\end{eqnarray}
Define the scattering transform by
\begin{equation}\label{scattering-transform}
\vct{t}(\lambda) = \int_{\C}e_{\lambda}(z)q_0(z)\mu(z,\lambda)\dee z_1\dee z_2
\end{equation}
and the operators
\begin{eqnarray}
\TT:\quad \TT f(z,\lambda) &=& \sgn(\abs{\lambda}^2-1)\frac{\vct{t}(\lambda)}{4\pi\kk{\lambda}}e_{-\lambda}(z)\kk{f(z,\lambda)}\label{T-def}\\
\CC:\quad \CC f(z,\lambda) &=& \frac{1}{\pi}\int_{\C}\frac{f(z,w)}{w-\lambda}\dee w_1\dee w_2\label{Cauchy}.
\end{eqnarray}
The following D-bar equation on the left holds, along with its integral form on the right:
\begin{equation}\label{D-bar}
\dbar_\lambda\mu(z,\lambda) = \TT\mu(z,\lambda),\qquad \mu(z,\lambda) = 1-\CC\TT\mu(z,\lambda).
\end{equation}

We still need some connections between the data $\Lambda_q$, the CGO solution $\mu(z,\lambda)$ and the potential $q_0(z)$. Define the operator
\begin{equation}\label{Slambda-def}
  (\SSS_\lambda \phi)(z) := \int_{\DOm}G_\lambda(z-y)\phi(y)\dee s(y),\qquad G_\lambda(z) = e^{\frac{\rmi\sqrt{E}}{2}(\lambda \kk{z}+\frac{z}{\lambda})}g_\lambda(z).
\end{equation}
The CGO solution $\psi(z,\lambda)$ satisfies the boundary integral equation
\begin{equation}\label{BIE-lambda}
(I+\SSS_\lambda(\Lambda_q-\Lambda_{-E}))\psi(\cdot,\lambda)|_{\DOm} = e^{\frac{\rmi\sqrt{E}}{2}(\lambda\kk{z}+\frac{1}{\lambda}z)}|_{\DOm},
\end{equation}
where the DN-map $\Lambda_{-E}$ corresponds to the potential $q_0 = 0$. In conjunction we have
\begin{equation}\label{scat-trans-DN}
\vct{t}(\lambda) = \int_{\bound}e^{\frac{-\rmi\sqrt{E}}{2}(\kk{\lambda}z+\kk{z}/\kk{\lambda})}(\Lambda_q-\Lambda_{-E})\psi(z,\lambda)\dee s(z).
\end{equation}
The potential can be reconstructed by
\begin{equation}\label{q0-recon1}
 q_0(z) = 2\rmi\sqrt{E}\partial_z\mu_{-1}^{\infty}(z),\qquad  \mu(z,\lambda) = 1+\frac{\mu_{-1}^{\infty}(z)}{\lambda}+\OO(\frac{1}{\abs{\lambda}}).
\end{equation}
Thus, in the absence of exceptional points, we have the necessary steps to reconstruct the potential $q_0$ from the DN-map $\Lambda_q$:
\begin{enumerate}
\item Solve $\psi(\cdot,\lambda)|_{\DOm}$ from the boundary integral equation \eqref{BIE-lambda}.
\item Compute the scattering transform using \eqref{scat-trans-DN}.
\item Choose a reconstruction point $z'$ and solve $\mu(z',\lambda)$ from the integral equation of \eqref{D-bar}.
\item Compute $q_0(z')$ from \eqref{q0-recon1}.
\end{enumerate}
This procedure fails if there are exceptional points. It has been known that when the potential is small compared to the energy, we have no exceptional points since the LS equation \eqref{LS-lambda} can be solved using the Neumann series. In the zero-energy case $E=0$ Nachman proved the absence of exceptional points when the potential is of \textit{conductivity type} \cite{Nachman1996}. In \cite{Music2013} this concept, at zero energy, was further studied, including a numerical test to find exceptional points for radial potentials perturbed from the conductivity-type. In \cite{deHoop2015} the same numerical test was carried out in the case of $E>0$. In our paper we again conduct the same test at $E<0$: the results of section \ref{sec:excep} are the first (numerical) results on exceptional points at negative energy, when the potential is not small. Also see \cite{Lakshtanov2015a} for similar results for non-radial potentials.

\section{Numerical implementation of the D-bar method}\label{sec:Dbarnumerics}
The following is the standard way of implementing the D-bar method and of conducting numerical tests, for $E=0$ see \cite{Siltanen2000,Siltanen2001} and for $E>0$  see \cite{deHoop2015}. Further, see the book by Mueller and Siltanen \cite{Mueller2012} which gives a detailed outlook. The basis of our method is similar in all three energy cases.

We use truncated Fourier basis on the boundary of the unit disk in order to approximate the operators $\Lambda_q$, $\Lambda_{-E}$ and $\SSS_\lambda$ by finite matrices $\mtx{L}_q$, $\mtx{L}_{-E}$ and $\mtx{S}_\lambda$ respectively. Choose an integer $N>0$ and define the basis functions
\begin{equation}\label{basisfunctions}
\phi^{(n)}(\theta) = \frac{1}{\sqrt{2\pi}}e^{\rmi n\theta},\quad n=-N,...,N.
\end{equation}
For the DN-map, solve the problem
\begin{equation}\label{gelfand}
  (-\Delta+q) u^{(n)} = 0\mbox{ in }\Omega,\qquad
  u^{(n)} = \phi^{(n)} \mbox{ on }\DOm
\end{equation}
for $u^{(n)}$ using Finite Element Method. Define $\mtx{L}_q=[\widehat{u}(\ell, n)]$ by
\begin{equation}\label{Lq-definition}
   \widehat{u}(\ell, n) = \int_{\DOm} \frac{\partial u^{(n)}}{\partial\nu} \overline{\phi^{(\ell)}}\dee s.
\end{equation}
Here $\ell$ is the row index and $n$ is the column index. The integration can be computed when the set $[0,2\pi)$ is divided into discrete points. The matrix $\mtx{L}_q$ represents the operator $\Lambda_q$ of \eqref{DNq-def} approximately. We add simulated measurement noise by defining
\begin{equation}\label{noiselevel}
  \mtx{L}_q^\epsilon := \mtx{L}_q + c\mtx{G},
\end{equation}
where $\mtx{G}$ is a $(2N+1)\times (2N+1)$ matrix with random entries
independently distributed according to the Gaussian normal density $\mathcal{N}(0,1)$. The constant $c>0$ can be adjusted for different relative errors $\|\mtx{L}_q^\epsilon-\mtx{L}_q\|/\norm{\mtx{L}_q}$, where $\norm{\cdot}$ is the standard matrix norm. The DN-map $\Lambda_{-E}$ is represented by the matrix $\mtx{L}_{-E}$ in a similar way, in the boundary value problem \eqref{gelfand} we then have $q=-E$. 

The matrix representation $\mtx{S}_\lambda$ of the operator $\SSS_\lambda$ of \eqref{Slambda-def} is obtained in a similar way. It is then possible to approximate the boundary integral equation \eqref{BIE-lambda} by the matrix equation
\begin{equation}\label{BIE-discrete}
(\textrm{I}+\mtx{S}_\lambda(\mtx{L}_q-\mtx{L}_{-E}))\vct{\psi}_\lambda = \vct{e}_\lambda,
\end{equation}
where $\textrm{I}$ is the correct sized unit matrix. In principle, this is solved for the vector $\vct{\psi}_\lambda$ by inverting the matrix $\textrm{I}+\mtx{S}_\lambda(\mtx{L}_q-\mtx{L}_{-E})$. We denote by $\mathcal{F}^{-1}$ the transformation from the Fourier series domain to the function domain and simply use \eqref{scat-trans-DN} to get $\vct{t}(\lambda)$:
\begin{equation}\label{scat_matrix}
\T(\lambda) = \int_{\bound}e^{\frac{-\rmi\sqrt{E}}{2}(\kk{\lambda}z+\kk{z}/\kk{\lambda})}\mathcal{F}^{-1}((\mtx{L}_q-\mtx{L}_{-E})\vct{\psi}_\lambda)\dee s(z).
\end{equation}
However, a truncation of $\T(\lambda)$ is needed, since large values of $\lambda$ result to the non-solvability of \eqref{BIE-discrete} because of exponential terms in the scattering transform \eqref{scat-trans-DN} and in the operator \eqref{Slambda-def}. We put the scattering transform to zero outside of an ellipse with the radius
\begin{equation}\label{ellipse}
r(\theta) = \frac{\sqrt{2}ab}{\sqrt{(b^2-a^2)\cos(2\theta-2\phi)+a^2+b^2}},
\end{equation}
where $a$ and $b$ are the semidiameters and the ellipse is rotated by $\phi$. An automatic choice of truncation of the scattering data is outside the scope of this initial feasibility study. Instead, the parameters of the ellipse are chosen in each case separately by looking at the scattering transform. 

We have the symmetry $\vct{t}(1/\kk{\lambda}) = \vct{t}(\lambda)$ which can be used to construct the scattering transform inside the unit circle. The truncated scattering transform used in the numerical simulations is then
\begin{equation}\label{scat-R}
\vct{t}_R(\lambda) = \left\{
\begin{array}{rcl}
                      0,&\quad           &\abs{\lambda}\leq 1/r(\theta) \\
\vct{t}(1/\kk{\lambda}),& \quad 1/r(\theta)\leq&\abs{\lambda}<1/R_1 \\
                      0,&\quad  1/R_1\leq&\abs{\lambda}\leq R_1 \\
       \vct{t}(\lambda),& \quad      R_1<&\abs{\lambda}< r(\theta)\\
                      0,&\quad           &\abs{\lambda}\geq r(\theta).
\end{array}
\right.
\end{equation}
The radius $R_1$ is needed since values of $\abs{\lambda}$ near the unit circle will also result to computational problems in the numerical Faddeev Green's function $g_\lambda(z)$. This essentially regularizes the D-bar method to handle noisy data, see the analogous case of zero-energy \cite{Knudsen2009}. 

After the truncation we transform to the modified version of the integral form of \eqref{D-bar},
\begin{equation}\label{IEIS-modified}
\mu_R = 1-\CC\TT_R\mu_R,
\end{equation}
where $\TT_R$ is the operator of \eqref{T-def} with $\T_R(\lambda)$ instead of $\T(\lambda)$. Equation \eqref{IEIS-modified} is solved by periodization and GMRES as explained in \cite{Knudsen2004}.

Choose a reconstruction point $z_r$. Let $dz$ be the finite difference and define the points 
$$
z_1 = z_r+dz,\quad z_2 = z_r-dz,\quad z_3 = z_r+\rmi\cdot dz,\quad z_4 = z_r-\rmi\cdot dz.
$$
Using the earlier described method we can solve the corresponding CGO solutions $\mu_R^i = \mu_R(z_i,\lambda)$, $i=1,2,3,4$. We omit the term $\OO(1/\abs{\lambda})$ in \eqref{q0-recon1}, use a finite $\lambda$ and finite difference method for the differentiation to get the approximate reconstruction equation
\begin{equation}\label{q0recon-discrete}
q_0(z_r) \approx \lambda\sqrt{E}\left(\frac{\mu_R^2-\mu_R^1}{2dz}+\rmi\frac{\mu_R^3-\mu_R^4}{2dz}\right).
\end{equation} 
Note that the CGO solutions are computed in a grid of parameters $\lambda$ and thus the same applies to the reconstruction $q_0(z_r)$. We can take the average value over indices corresponding to different values of $\lambda$.

\subsection{Error caused by the truncation}
The following theoretical results readily follow from \cite{Santacesaria2013}. In the case $E>0$ in \cite{deHoop2015} similar results were proved, but here we write them out more clearly. In \cite{Santacesaria2013} the stability of the D-bar method at negative energy, essentially of logarithmic type, was rigorously analyzed assuming the following additional properties to us:
\begin{equation}\label{assumptions}
\begin{split}
&q_0\in W^{m,1}(\R^2)\textrm{ for some }m>2,\qquad \abs{E}>E_1 = E_1(\norm{q_0}_{m,1},\Omega),\\
&\norm{q_0}_{m,1}\leq N \quad \textrm{ where }\norm{f}_{m,1} = \max_{\abs{J}\leq m}\norm{\partial^{J}f}_{L^1(\R^2)}.
\end{split}
\end{equation}
Let $R$ be a truncation parameter large enough such that 
\begin{equation*}
R>\max(1,2R'/\sqrt{\abs{E}}),\qquad R<\min(1,\sqrt{\abs{E}}/(2R')),
\end{equation*}
where $R'$ is from Lemma 2.1 of \cite{Santacesaria2013}. Consider the simplified truncated scattering transform $\vct{t}_{R,s}(\lambda)$, where $\vct{t}_{R,s}(\lambda) = 0$ for $\abs{\lambda}<1/R$ and for $\abs{\lambda}>R$. Lemma 3.1 of \cite{Santacesaria2013} gives us
\begin{equation}\label{L3.1}
\begin{split}
\norm{\abs{\lambda}^j\frac{\vct{t}(\lambda)}{\kk{\lambda}}}_{L^p(\abs{\lambda}<1/R)} &\leq C(N,m,p)\abs{E}^{-m/2}R^{-(m-1+j+2/p)},\\
\norm{\abs{\lambda}^j\frac{\vct{t}(\lambda)}{\kk{\lambda}}}_{L^p(\abs{\lambda}>R)} &\leq C(N,m,p)\abs{E}^{-m/2}R^{-m-1+j+2/p},
\end{split}
\end{equation}
where $C(N,m,p)$ is a constant depending on $N,m$ and $p$ and $j=-1,0,1$. Denote 
$$
\A_R = \{\lambda\in\C:\quad 1/R\leq\abs{\lambda}\leq R\}.
$$
\begin{lemma}\label{apulemma}
Let $m>2$ and $p\geq 1$ as assumed in \eqref{assumptions}. Then for $j=-1,0,1$ we have
\begin{equation}
\norm{\abs{\lambda}^j\frac{\vct{t}(\lambda)}{\kk{\lambda}}}_{L^p(\C\setminus \A_R)} \leq C(N,m,p)\abs{E}^{-m/2}R^{-m+2}.
\end{equation}
\begin{proof}
We have $0<1-j<2$, $-2<-2/p<0$, $-2<-1+j<0$ and $0<2/p<2$. Using these and \eqref{L3.1} we get
\begin{equation}
\begin{split}
\norm{\abs{\lambda}^j\frac{\vct{t}(\lambda)}{\kk{\lambda}}}_{L^p(\C\setminus \A_R)} &\leq C(N,m,p)\abs{E}^{-m/2}(R^{-m+1-j-2/p)}+R^{-m-1+j+2/p})\\
& \leq C(N,m,p)\abs{E}^{-m/2}(R^{-m-2}+R^{-m+2})\\
& \leq C(N,m,p)\abs{E}^{-m/2}\cdot R^{-m+2}.
\end{split}
\end{equation}
\end{proof}
\end{lemma}
Note that in this subsection we used a simpler truncation than \eqref{scat-R}, and so the result below would have to be fine-tuned to be exact. However, the idea of truncation is the same in both versions of truncation, and so the use of \eqref{scat-R} is justified by the following theorem.
\begin{theorem}
Assume \eqref{assumptions}. Let $q_R$ be the reconstructed potential using $\vct{t}_{R,s}$. We have
\begin{equation}\label{eq:q0qR}
\abs{q_0(z)-q_R(z)}\leq C(\Omega,N,m,p)\abs{E}^{(-m+2)/2}R^{-m+2}.
\end{equation}
\begin{proof}
See \cite{Santacesaria2013} last estimate before (4.8), using our notation it reads
\begin{equation}
\begin{split}
\abs{q(z)-q_R(z)}\leq & C(\Omega,N,m,p)\abs{E}\biggl(\norm{\abs{\lambda}\frac{\vct{t}(\lambda)}{\kk{\lambda}}}_{L^1(\C\setminus \A_R)}+\norm{\frac{1}{\abs{\lambda}}\frac{\vct{t}(\lambda)}{\kk{\lambda}}}_{L^1(\C\setminus \A_R)} \\
& + \sum_{j=-1}^1\norm{\abs{\lambda}^j\frac{\vct{t}(\lambda)}{\kk{\lambda}}}_{L^p(\C\setminus \A_R)}+ \abs{E}^{-1/2}\norm{\frac{\vct{t}(\lambda)}{\kk{\lambda}}}_{L^p(\C\setminus \A_R)}\biggl).
\end{split}
\end{equation}
We take the term $\abs{E}^{-1/2}$ out by multiplying by a necessary constant and use Lemma \ref{apulemma} to all of the terms which gives \eqref{eq:q0qR}.
\end{proof}
\end{theorem}

\subsection{Computation of the Faddeev Green's function at negative energy}\label{sec:green}
In order to numerically compute the boundary integral equation \eqref{BIE-discrete}, or to compute the CGO solution directly from \eqref{LS-lambda}, we need a numerical algorithm for $g_\lambda(z)$ for any $\lambda\in\C\setminus D(0,1)$, $z\in D(0,1)$ and of course $E<0$. The algorithm to follow stems originally from the zero-energy case and the numerical implementation presented in \cite{Siltanen1999}. For more details we refer to the positive-energy case \cite{deHoop2015}.

Following \cite{deHoop2015} we have, using $\zeta$-parameters,
\begin{equation}\label{gz1}
g_\zeta(z) = \frac{1}{2\pi}e^{-\rmi x_1k_1}\re(\int_{0}^{\infty}e^{\rmi x_2t}\frac{e^{-x_1\sqrt{(t+k_2\rmi)^2-E}}}{\sqrt{(t+k_2\rmi)^2-E}}\dee t),
\end{equation}
where $x_1\geq0$, $\zeta$ is in the reduced form $\zeta = \matr{k_1,0}^{\mathrm{T}}+\rmi\matr{0,k_2}^{\mathrm{T}}$ and $k_2>\abs{k_1}>0$. The integrand converges quickly for large $x_1$ which is a good feature, but it oscillates the larger $\abs{x_2}$ is, which is an unwanted feature. By complexifying and choosing suitable integration paths we can prove the following lemmata.
\begin{lemma}
Let $x_1\geq0$, $\zeta$ in the reduced form and $x_2\geq 0$. Then
\begin{equation}\label{gz2}
g_\zeta(z) = \frac{1}{2\pi}e^{-\rmi x_1k_1}\re(\int_0^\infty e^{-t}\frac{e^{-x_1 \rmi\sqrt{t^2/x_2^2+tk_2/x_2+k_1^2}}}{\sqrt{t^2/x_2^2+tk_2/x_2+k_1^2}}\dee t),
\end{equation}
\end{lemma}
\begin{lemma}
Let $x_1\geq0$, $\zeta$ in the reduced form and $x_2<0$. Then
\begin{equation}\label{gz3}
g_\zeta(z) = \frac{1}{2\pi}e^{-\rmi x_1k_1}\re(I_1-\rmi e^{\rmi x_2}\int_0^\infty e^{x_2 t}\frac{e^{-x_1\sqrt{(1+(k_2-t)\rmi)^2-E}}}{\sqrt{(1+(k_2-t)\rmi)^2-E}}\dee t),
\end{equation}
where
$$
I_1 = \int_{0}^{1}e^{\rmi tx_2}\frac{e^{-x_1\sqrt{(t+k_2\rmi)^2-E}}}{\sqrt{(t+k_2\rmi)^2-E}}.
$$
\end{lemma}
The integrands in \eqref{gz2} and \eqref{gz3} converge quickly for large $\abs{x_2}$ complementing the formula \eqref{gz1}. Note however that for any of these formulas a small $z=x_1+\rmi x_2$ will be a problem because of slow convergence.

Write $g_\zeta^{T_1}$,$g_\zeta^{T_2}$ and $g_\zeta^{T_3}$ for the finite integrals \eqref{gz1},\eqref{gz2} and \eqref{gz3} respectively. We need to choose the upper limits $T_i$, $i=1,2,3$. For the three cases we find that the remainder of the integral is limited by the exponential integral
$$
Ei(t) = \int_{t}^{\infty}\frac{e^{-s}}{s}\dee s.
$$
Then, similarly to \cite{deHoop2015}, we can obtain $ \abs{g_\zeta-g_\zeta^{T_i}} < 1E-8$ by choosing
\begin{eqnarray}
T_1 &=& \max\{\frac{14\sqrt{2}}{x_1c_1},\sqrt{2}\abs{k_1}\}\label{T1-def},\\
T_2 &=& 14\label{T2-def},\\
T_3 &=& \frac{14}{c_2x_1-x_2}\label{T3-def},
\end{eqnarray}
where
\begin{eqnarray*}
c_1 &=& \cos(\theta_1),\quad \theta_1 = Arg(\sqrt{k_1^2+2\sqrt{2}\abs{k_1}k_2 \rmi}),\\
c_2 &=& \cos(\theta_2),\quad \theta_2 = Arg(\sqrt{1-k_1^2+2k_2 \rmi}).
\end{eqnarray*}
For $g_\zeta^{T_i}$, the integration range $[0,T_i]$ is divided into $M_i$ points (with $g_\zeta^{T_3}$ there is also the additional integral $I_1$) and the Gaussian quadrature is used. The integers $M_i$ are chosen large enough so that for any integer $M>M_i$ the first 8 digits are not changing in the numerical value of $g_\zeta^{T_i}(z)$.

\begin{figure}
\unitlength=1mm
\begin{picture}(120,60)
\epsfxsize=6cm
\put(0,0){\epsffile{\imagepath 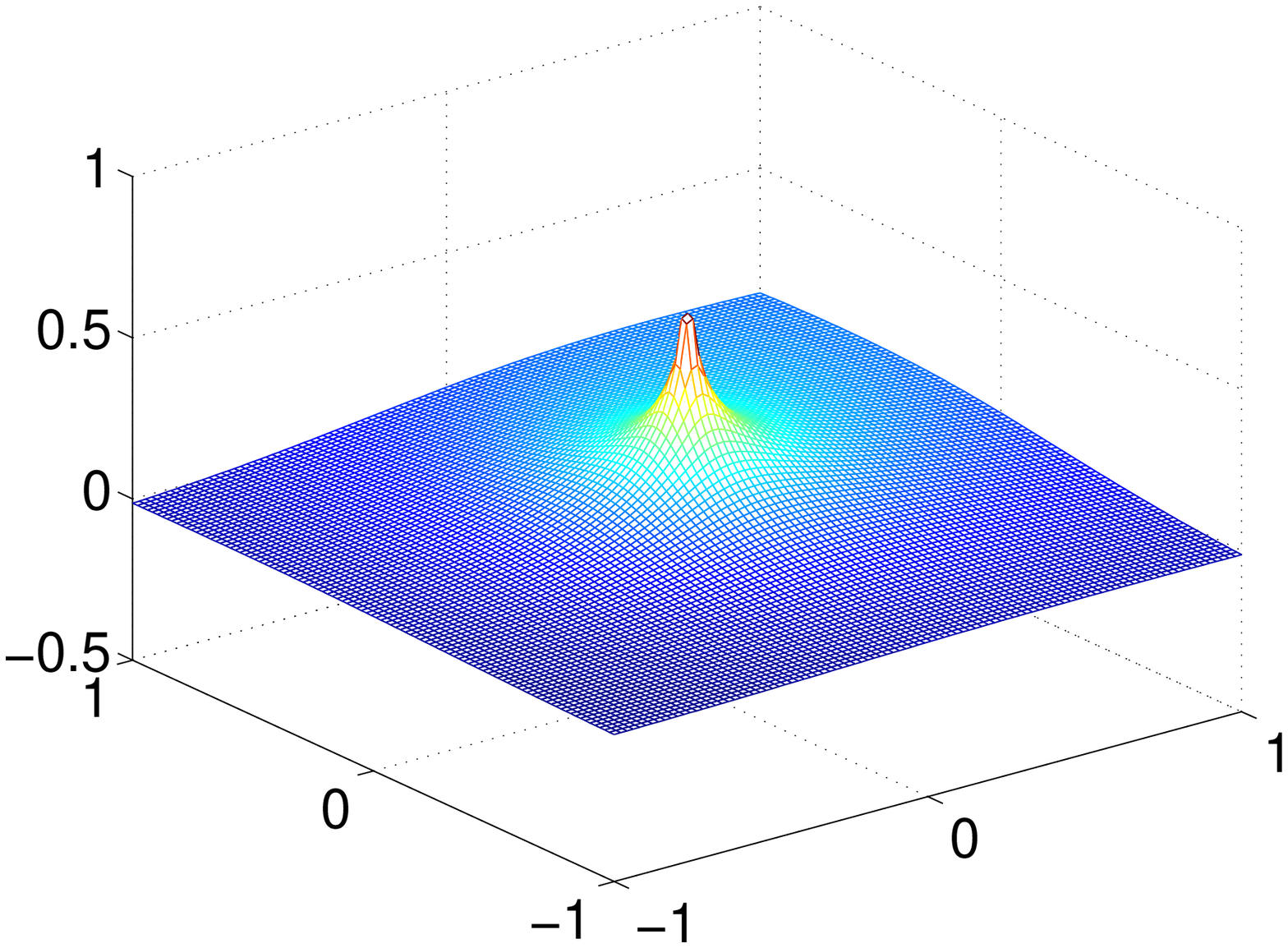}}
\epsfxsize=6cm
\put(62,0){\epsffile{\imagepath 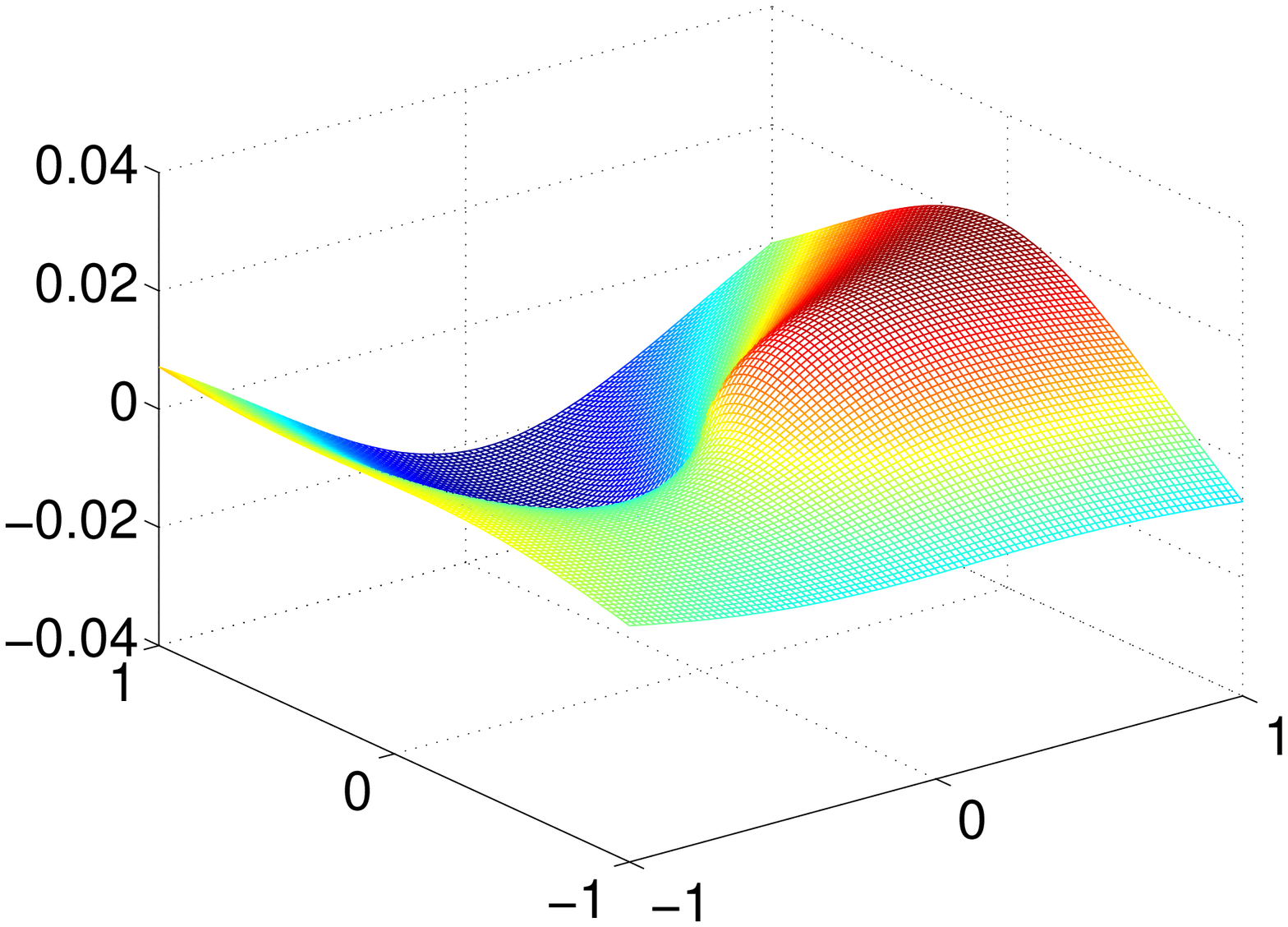}}
\put(90,49){\small $\im(g_\lambda(z))$}
\put(28,49){\small $\re(g_\lambda(z))$}
\end{picture}
\caption{\label{fig:gzeta}The real and imaginary parts of $g_{\lambda}(z)$ in $400\times400$ grid of points $z$, $\lambda=1+i$, $E=-1$.}
\end{figure}

See \cite{deHoop2015} for a scaling and switching relations of the Faddeev Green's function. The $z$-plane is divided into computational areas:
\begin{itemize}
\item For $\abs{z}<0.01$ we put $g_\lambda(z) = 0$.
\item For $0.01\leq \abs{z} <0.5$ we scale $z$ outwards from the origin by a factor of $100$.
\item For $0.5\leq \abs{z} <1$ we scale $z$ outwards from the origin by a factor of $2$.
\item For $1\leq \abs{z}$, $x_1<0$ we use the switching relation to obtain $x_1>0$.
\item For $1\leq \abs{z}$, $x_1\geq0$, $-x_1<x_2<0.5 x_1$ we use $g_\lambda^{T_1}(z)$ originating from \eqref{gz1}.
\item For $1\leq \abs{z}$, $x_1\geq0$, $x_2\geq 0.5 x_1$ we use $g_\lambda^{T_2}(z)$ originating from \eqref{gz2}.
\item For $1\leq \abs{z}$, $x_1\geq0$, $x_2<-x_1$ we use $g_\lambda^{T_3}(z)$ originating from \eqref{gz3}.
\end{itemize}
These domains have been decided by numerical tests. A sample of the function $g_\lambda(z)$ is pictured in figure \ref{fig:gzeta} in $400\times400$ -grid of points $z$, $\lambda=1+i$, $E=-1$. This picture itself of course does not guarantee that our numerical implementation is correct. However it looks smooth meaning that the different computational areas are indistinguishable. In the next section a proper validation of the algorithm is conducted.

\subsection{Validation of the numerical Faddeev Green's function}
\begin{figure}
\unitlength=1mm
\begin{picture}(120,65)
\epsfxsize=5cm
\put(0,0){\epsffile{\imagepath 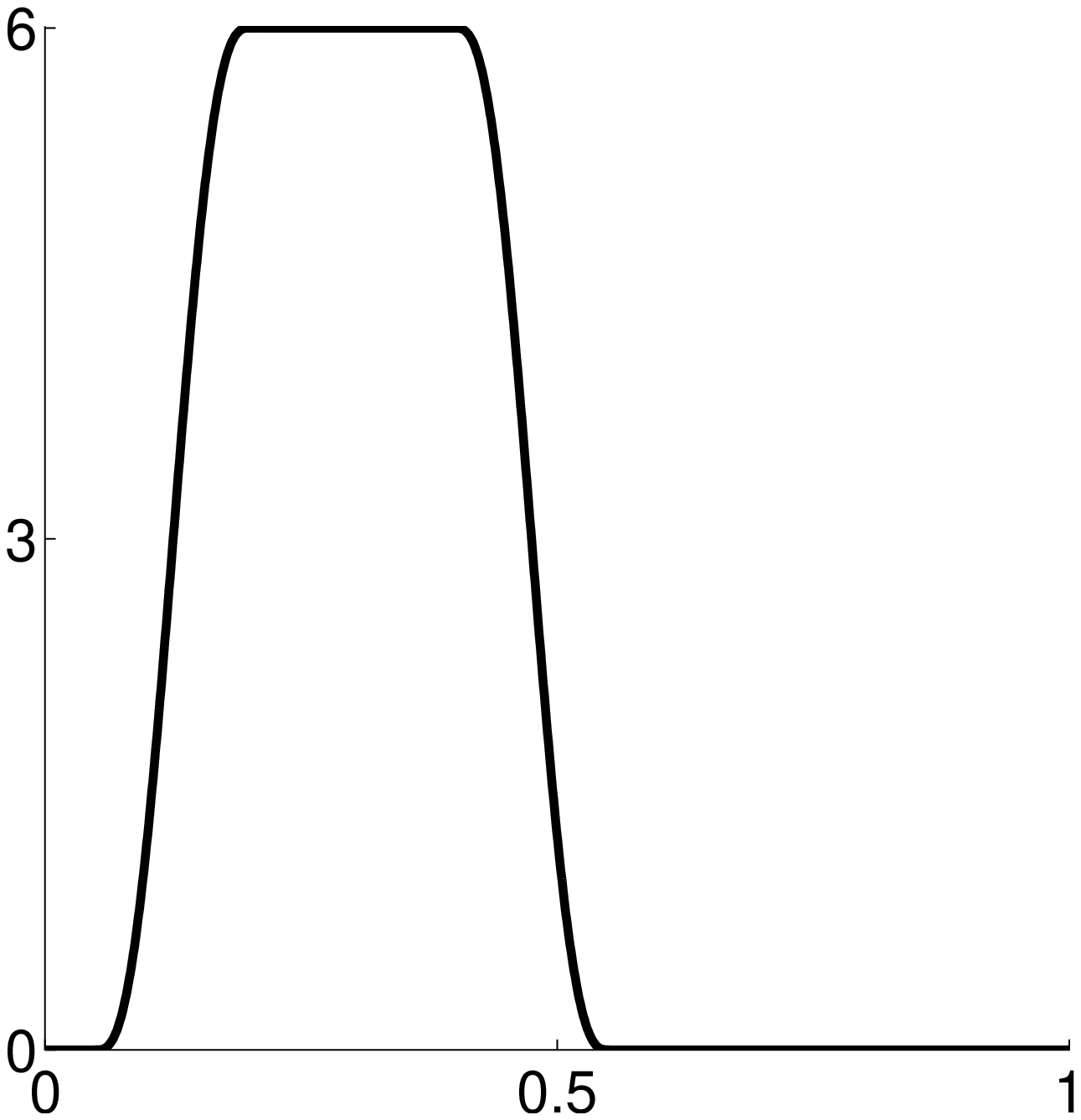}}
\epsfxsize=6cm
\put(60,0){\epsffile{\imagepath 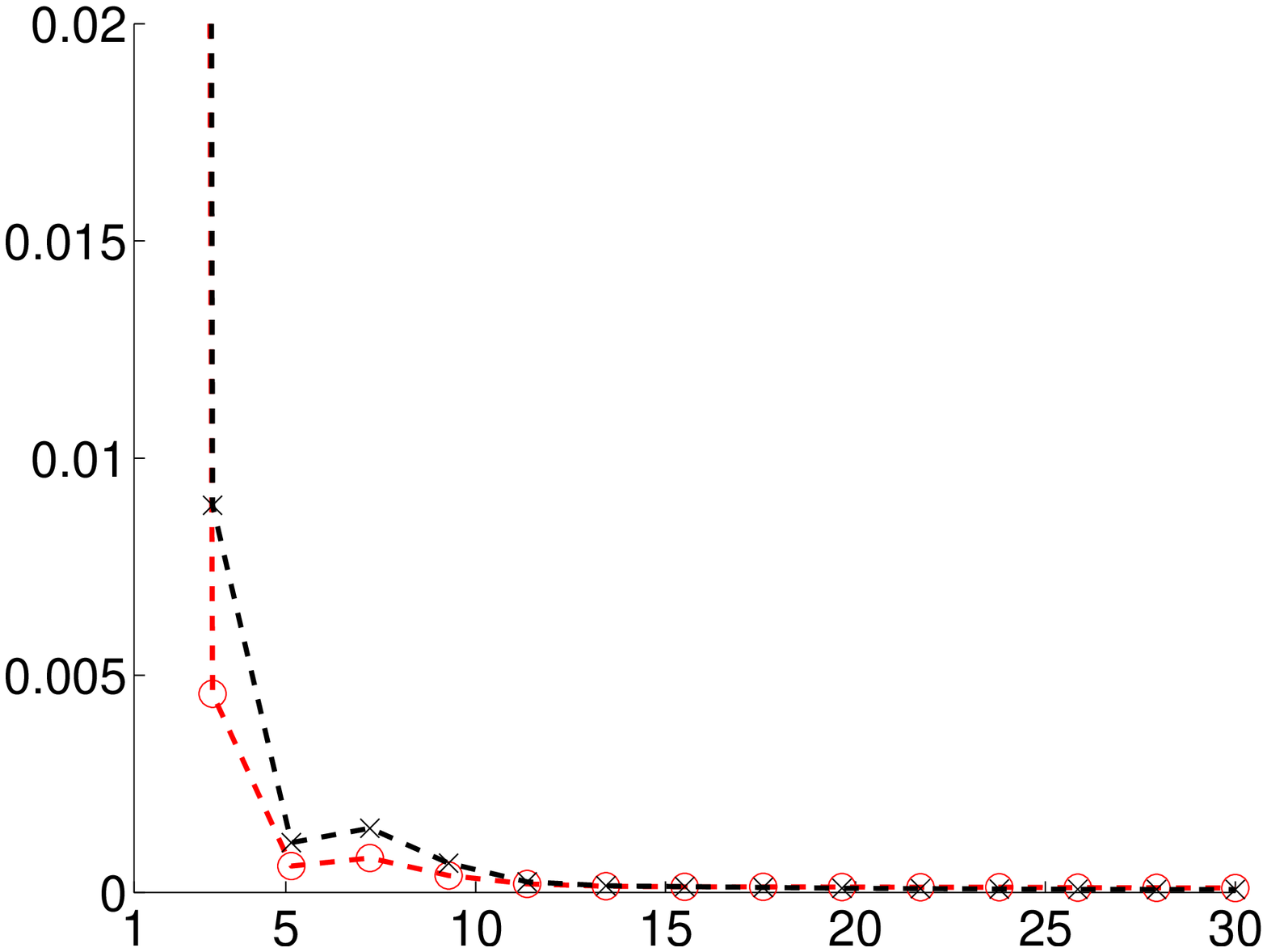}}
\put(0,52){\small $q_0(\abs{z})$}
\put(50,-3){\small $|z|$}
\put(62,52){\tiny $\norm{\dbar \mu-\frac{1}{4\pi\kk{\lambda}}\vct{t}(\lambda)e_{-\lambda}(z)\kk{\mu_0}}_{L^2(D(0,1))}$}
\put(117,-3){\small $|\lambda|$}
\end{picture}
\caption{\label{fig:greentest}On the left: the profile of a test potential used in the validation of the Green's function. On the right: the relative error \eqref{dbar-error} with $M=6$ in black and $M=7$ in red.}
\end{figure}

See the radial potential whose profile is pictured in figure \ref{fig:greentest} on the left. For this potential, we can solve the CGO solution $\mu$ from the LS equation \eqref{LS-lambda} and test the $\dbar$ -equation \eqref{D-bar} using the five-point stencil method. Let us use the finite difference of $d\lambda=0.0001$, parameters $\lambda$ from 1 to 30 and a $2^M\times 2^M$ grid in the $z$-space with $M=6$ and $M=7$ for the LS solver; the solver is based on the ideas of Vainikko \cite{Vainikko2000} and it is detailed in \cite{Knudsen2004}. For each $\lambda=\lambda_1+\lambda_2 \rmi$, we compute
\begin{enumerate}
\item The CGO solution $\mu_0$ in the z-grid, corresponding to the parameter $\lambda$.
\item The CGO solutions $\mu_1,\mu_2,\mu_3,\mu_4,\mu_5,\mu_6,\mu_7$ and $\mu_8$ using $\lambda+d\lambda$, $\lambda+2d\lambda$, $\lambda-d\lambda$, $\lambda-2d\lambda$, $\lambda+d\lambda \rmi$, $\lambda+2d\lambda \rmi$, $\lambda-d\lambda \rmi$ and $\lambda-2d\lambda \rmi$ respectively.
\item The radial scattering transform $\vct{t}(\abs{\lambda})$ with $\mu_0$ using \eqref{scattering-transform}.
\item The derivatives and the $\dbar$ -operation by
\begin{eqnarray*}
\partial_{\lambda_1}\mu &=& \frac{-\mu_2+8\mu_1-8\mu_3+\mu_4}{12d\lambda}\\
\partial_{\lambda_2}\mu &=& \frac{-\mu_6+8\mu_5-8\mu_7+\mu_8}{12d\lambda}\\
\dbar \mu &=& \frac{1}{2}(\partial_{\lambda_1}+\rmi\partial_{\lambda_2})\mu.
\end{eqnarray*}
\item The error
  \begin{equation}\label{dbar-error}
    \norm{\dbar \mu-\frac{1}{4\pi\kk{\lambda}}\vct{t}(\lambda)e_{-\lambda}(z)\kk{\mu_0}}_{L^2(D(0,1))}.
  \end{equation}
\end{enumerate}
In figure \ref{fig:greentest} on the right we see the error \eqref{dbar-error} as a function of $\lambda$. The error grows very large near $\abs{\lambda}=1$ which is an unwanted feature of our numerical method for $g_\lambda(z)$. Overall, the size of the error is larger than in the positive energy case which can be partially explained by the lack of the \textit{layer potential method} for small $z$, see \cite{deHoop2015} and note that we put $g_\lambda(z)=0$ for $\abs{z}<0.01$. Otherwise, we conclude that our numerical method for $g_\lambda(z)$ is valid.

\subsection{Reconstruction of potentials at negative energy}\label{sec:q0recon}
Let us fix $E=-1$. In the numerical tests that follow, for the DN- and $\mtx{S}_\lambda$ -matrices we use $N=16$, see the definition \eqref{Lq-definition}. We add gaussian noise to each element using equation \eqref{noiselevel} so that the relative matrix norm between the original DN-matrix and the noisy DN-matrix is 0.005\%. In the mesh for the FEM we have 1048576 triangles.

\begin{figure}
\unitlength=1mm
\begin{picture}(120,80)
\epsfxsize=6cm
\put(0,0){\epsffile{\imagepath v1_q11D.eps}}
\epsfxsize=6cm
\put(60,0){\epsffile{\imagepath 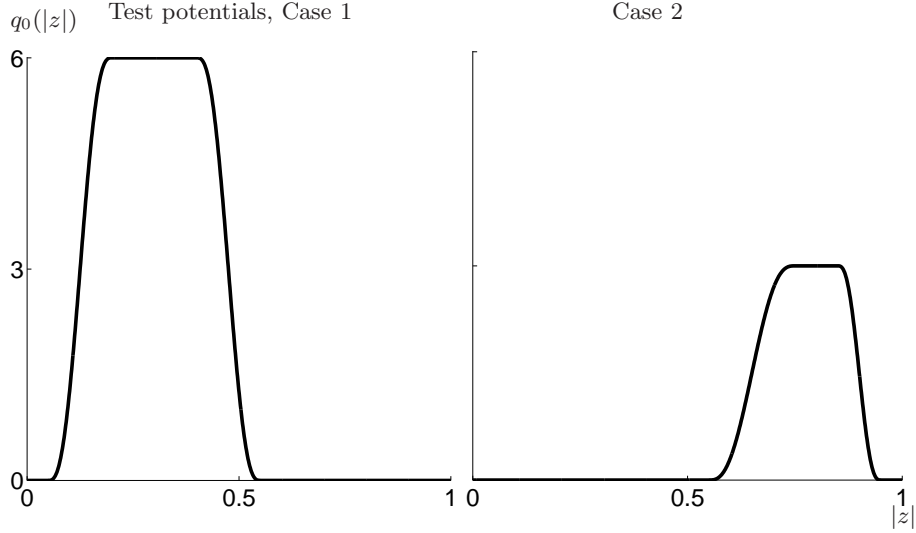}}
\put(13,65){\small Test potentials, Case 1}
\put(80,65){\small Case 2}
\put(0,64){\small $q_0(\abs{z})$}
\put(117,-2){\small $|z|$}
\end{picture}
\caption{\label{fig:q-profiles}Profile plots of Case 1  and Case 2 radial potentials $q_0(z)=q_0(\abs{z})$.}
\end{figure}

\begin{figure}[h]
\begin{picture}(120,60)
\epsfxsize=6cm
\put(0,1){\epsffile{\imagepath 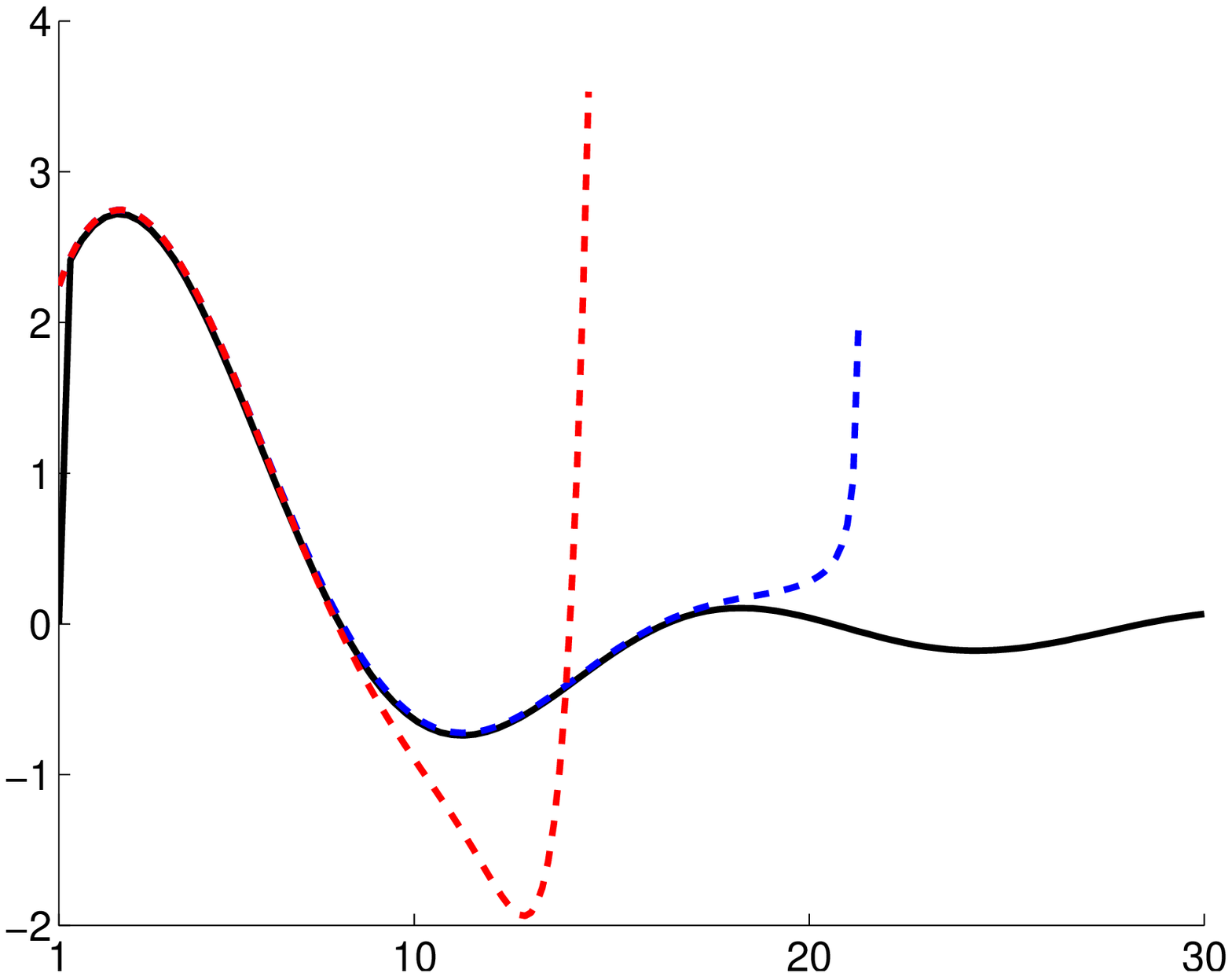}}
\epsfxsize=6cm
\put(60,1){\epsffile{\imagepath 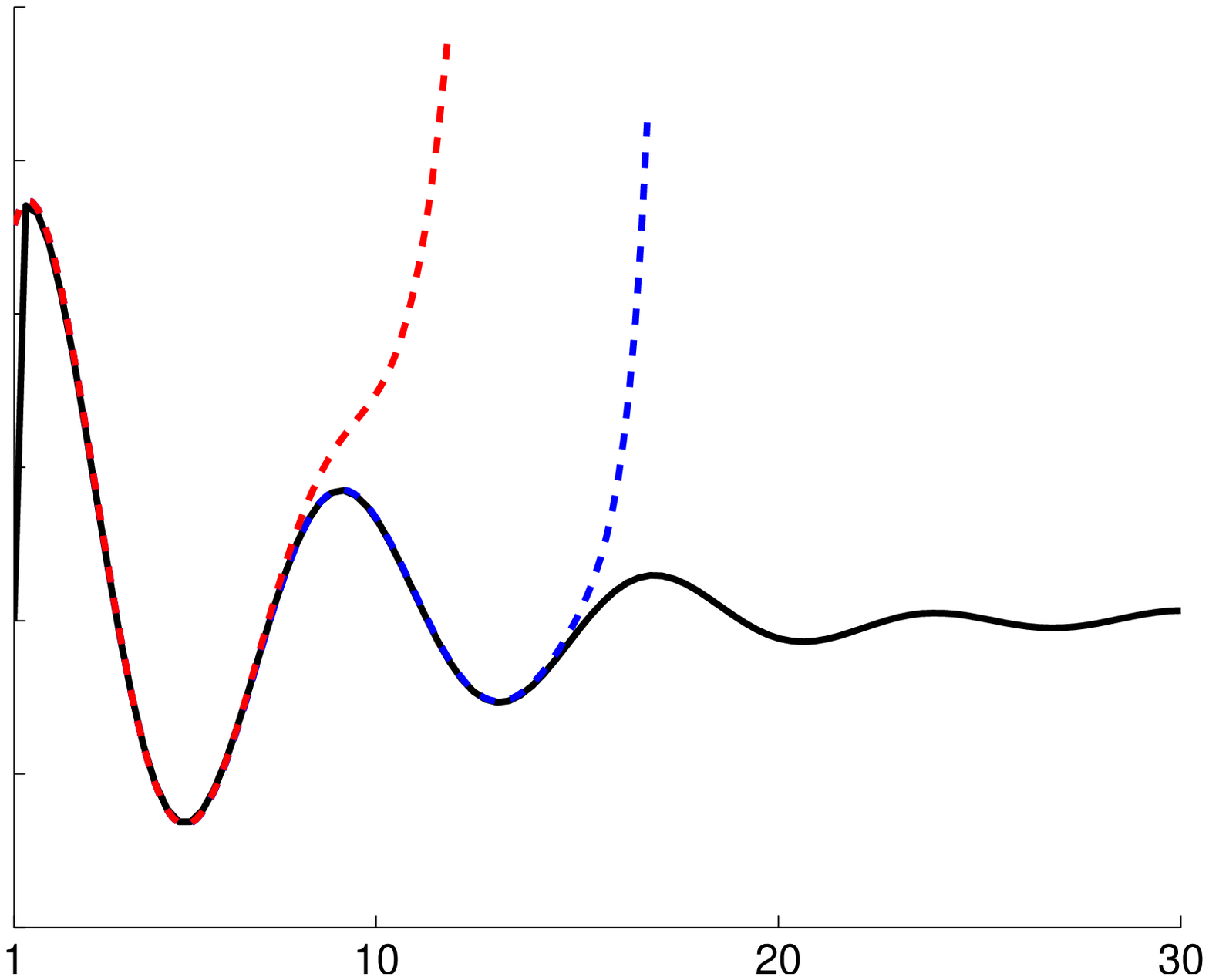}}
\put(26,53){Case 1}
\put(90,53){Case 2}
\put(120,-2){\small $\abs{\lambda}$}
\put(-1,53){\small $\vct{t}(\abs{\lambda})$}
\end{picture}
\caption{\label{scats}Used potentials are Case 1 and Case 2 of figure \ref{fig:q-profiles}, energy $E=-1$. The radial scattering transform $\vct{t}(\abs{\lambda})$ computed in three ways; black solid line indicates computation directly from \eqref{LS-lambda} using the knowledge of $q_0$, blue dashed line indicates computation using the DN-matrices $\mtx{L}_q$, $\mtx{L}_{-E}$ and using equations \eqref{BIE-discrete} and \eqref{scat_matrix}, red dashed line indicates the same but with noisy DN-matrix $\mtx{L}^\epsilon_q$.}
\end{figure}

We reconstruct two radially symmetric potentials (Case 1 and Case 2) of figure \ref{fig:q-profiles}. For radially symmetric potentials we have $\vct{t}(\lambda) = \vct{t}(\abs{\lambda})$ and $\im(\vct{t}(\lambda))=0$. In figure \ref{scats} we have the radially symmetric scattering transform computed in three ways; black solid line indicates computation directly from \eqref{LS-lambda} using the knowledge of $q_0$, blue dashed line indicates computation using the DN-matrices $\mtx{L}_q$ and $\mtx{L}_{-E}$ and using \eqref{BIE-discrete} and \eqref{scat_matrix}, red dashed line indicates the same but with noisy DN -matrix $\mtx{L}_q^\epsilon$. On the left we used Case 1 potential and on the right we used Case 2 potential.

\begin{figure}[h]
\begin{picture}(120,62)
\epsfxsize=6cm
\put(0,1){\epsffile{\imagepath 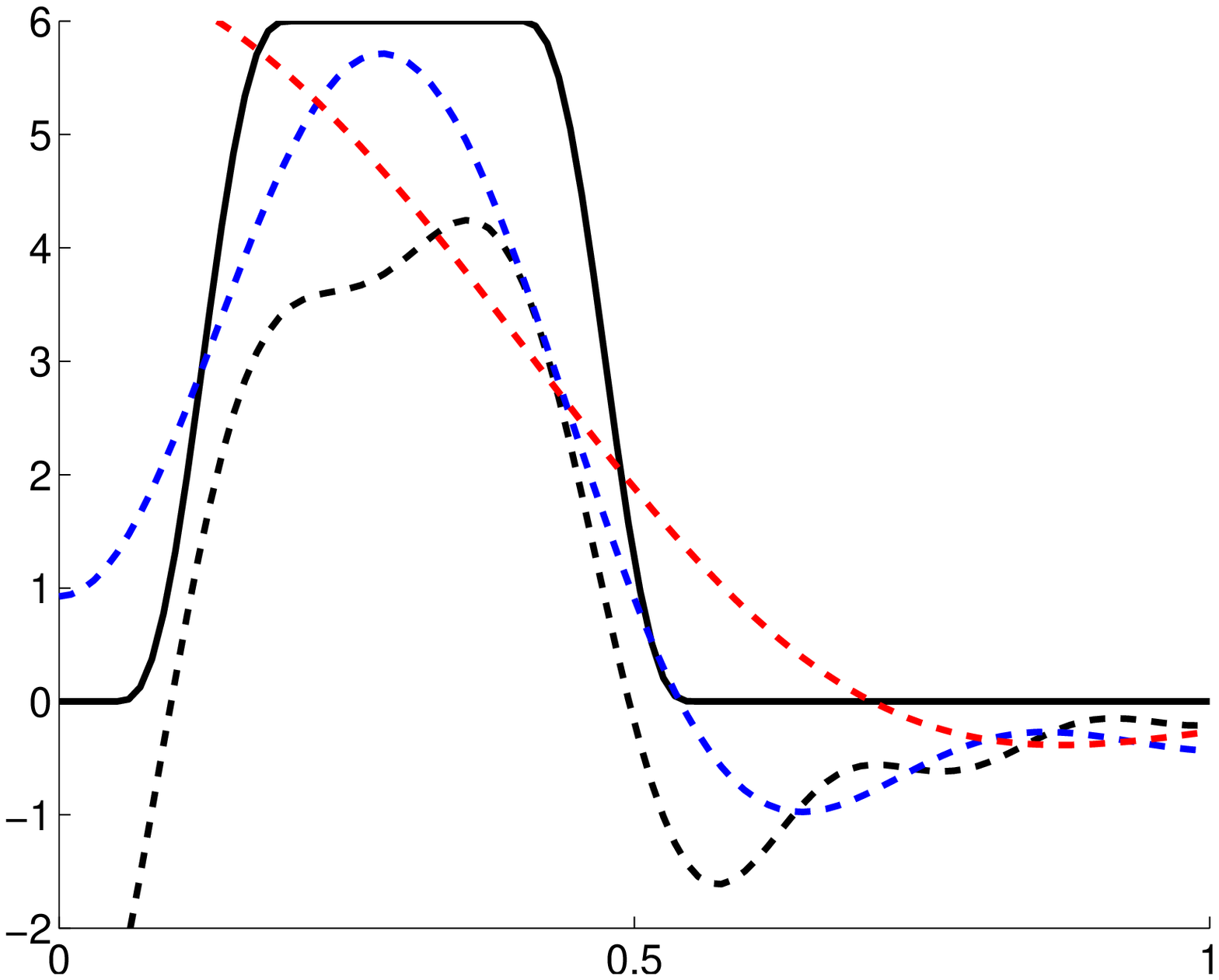}}
\epsfxsize=6cm
\put(60,1){\epsffile{\imagepath 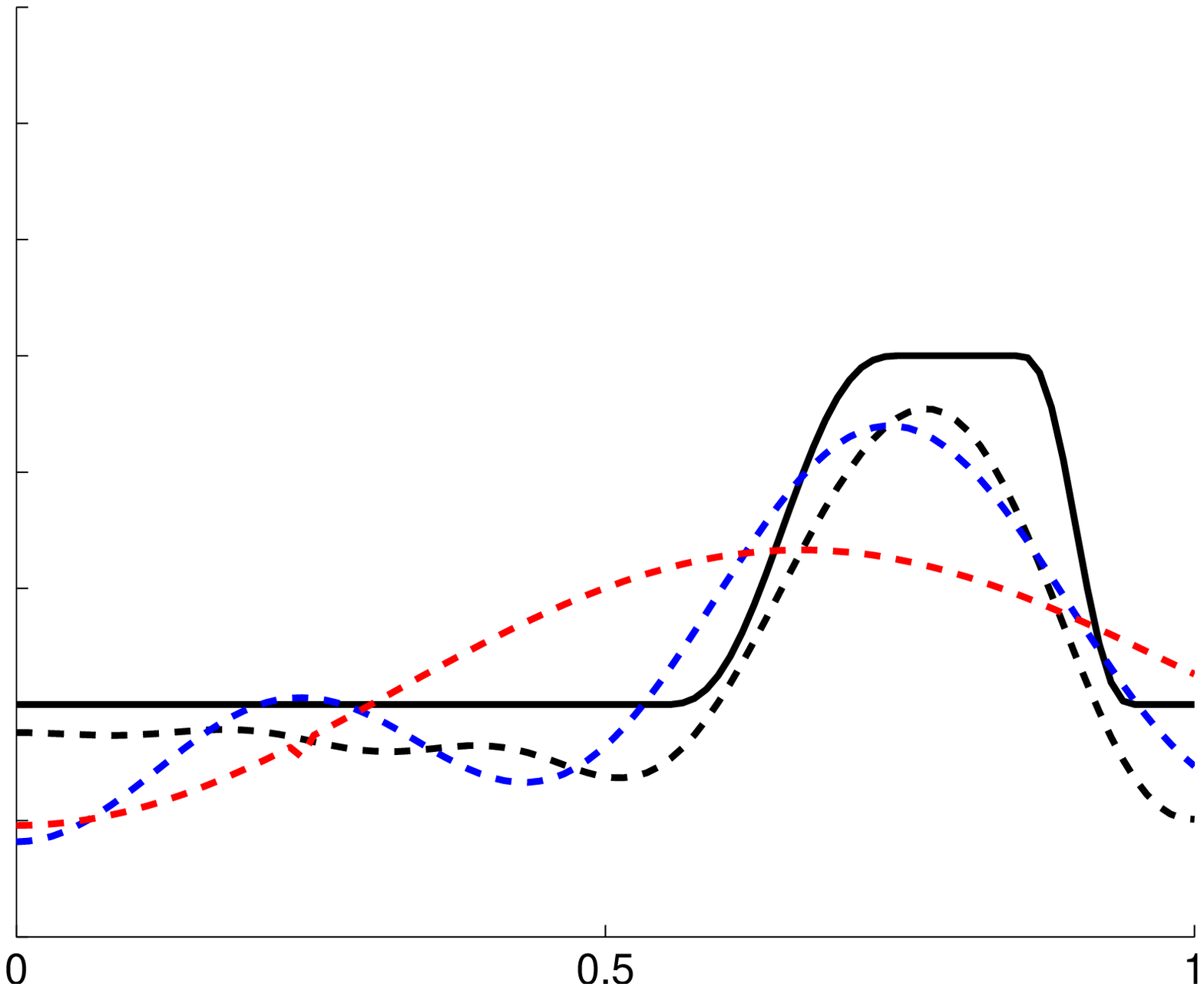}}
\put(26,53){Case 1}
\put(86,53){Case 2}
\put(-1,53){\small $q_0(\abs{z})$}
\put(120,-2){\small $\abs{z}$}
\end{picture}
\caption{\label{recon}Reconstruction of the radial potentials Case 1 and Case 2, energy $E=-1$. Black solid line: the original potential. Black dashed line: reconstruction using the knowledge of the potential via the scattering transform \eqref{scattering-transform}. Blue dashed line: reconstruction using the new computational method without noise. Red dashed line: reconstruction using the new computational method with added noise.}
\end{figure}

In figure \ref{recon} we plot the reconstructions of Cases 1 and 2 using the three different scattering transforms. Black solid line indicates the original potential, black dashed line indicates the reconstruction using the knowledge of the potential via the scattering transform \eqref{scattering-transform}, blue dashed line indicates reconstruction using \eqref{q0recon-discrete} without noise, and red dashed line indicates reconstruction using \eqref{q0recon-discrete} with added noise.

\section{Numerical investigation of exceptional points at negative energy}\label{sec:excep}
\begin{figure}
\unitlength=1mm
\begin{picture}(120,50)
\put(0,0){\includegraphics[width=4cm]{\imagepath 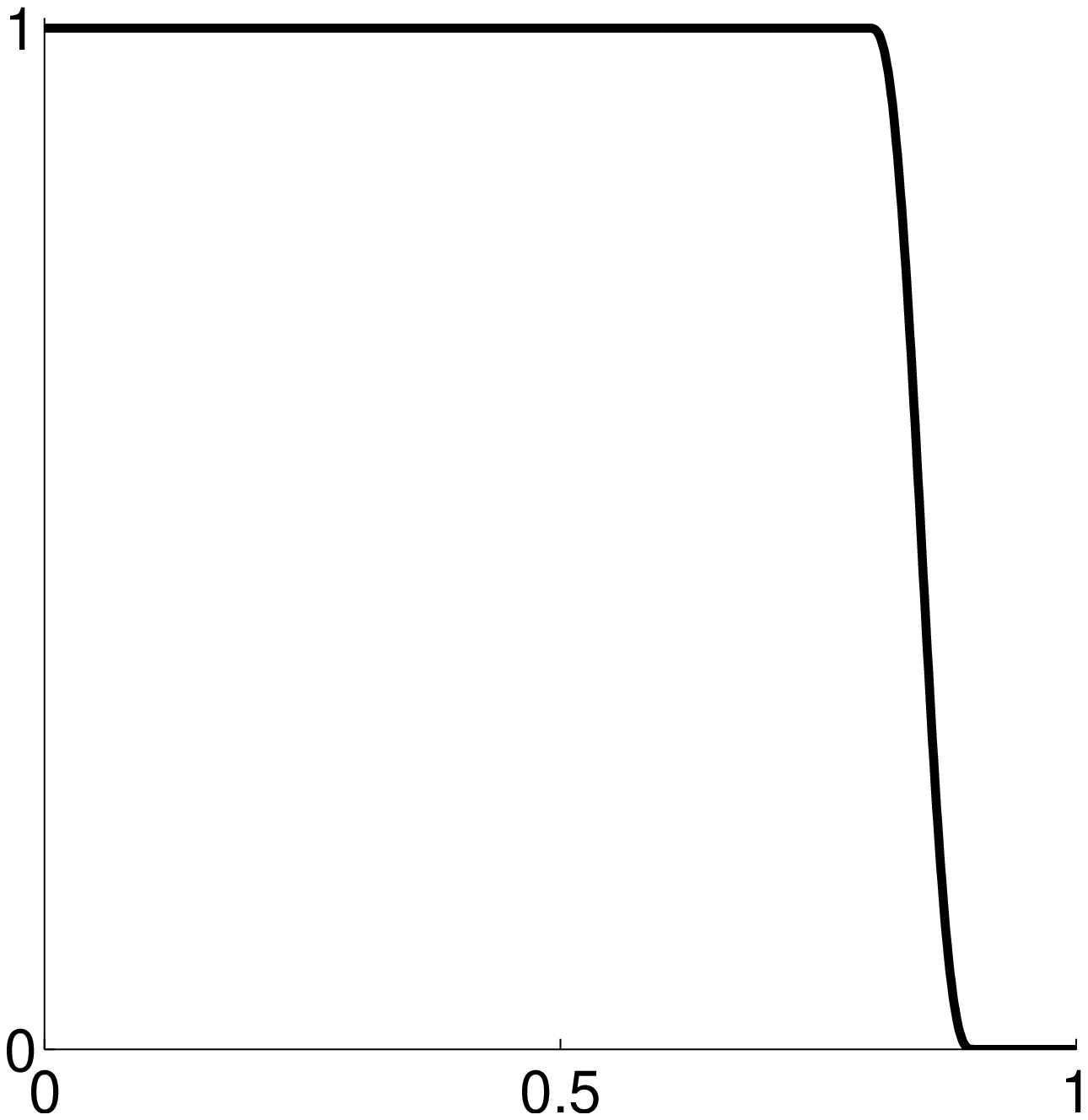}}
\put(40,0){\includegraphics[width=4cm]{\imagepath 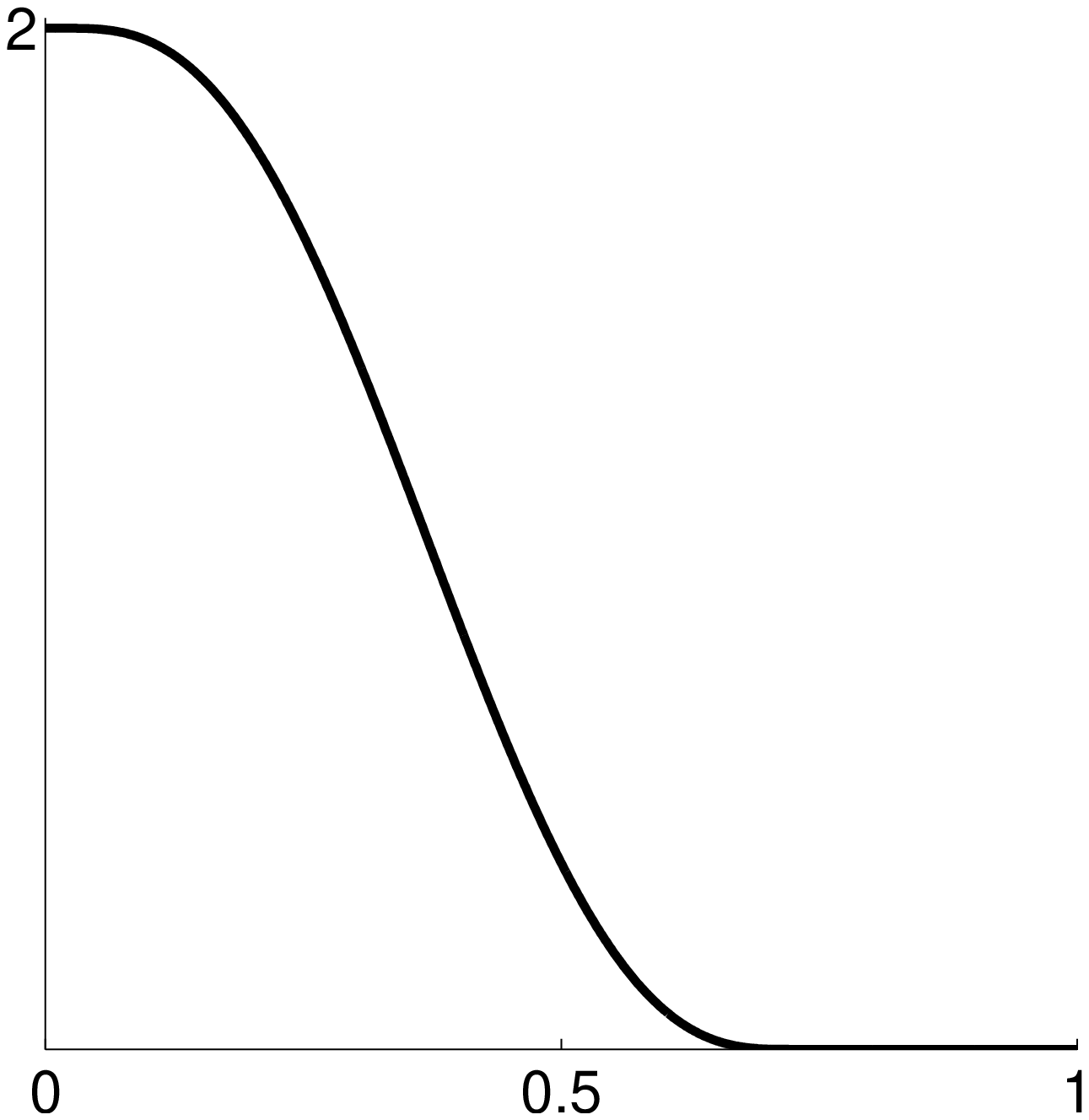}}
\put(80,0){\includegraphics[width=4cm]{\imagepath 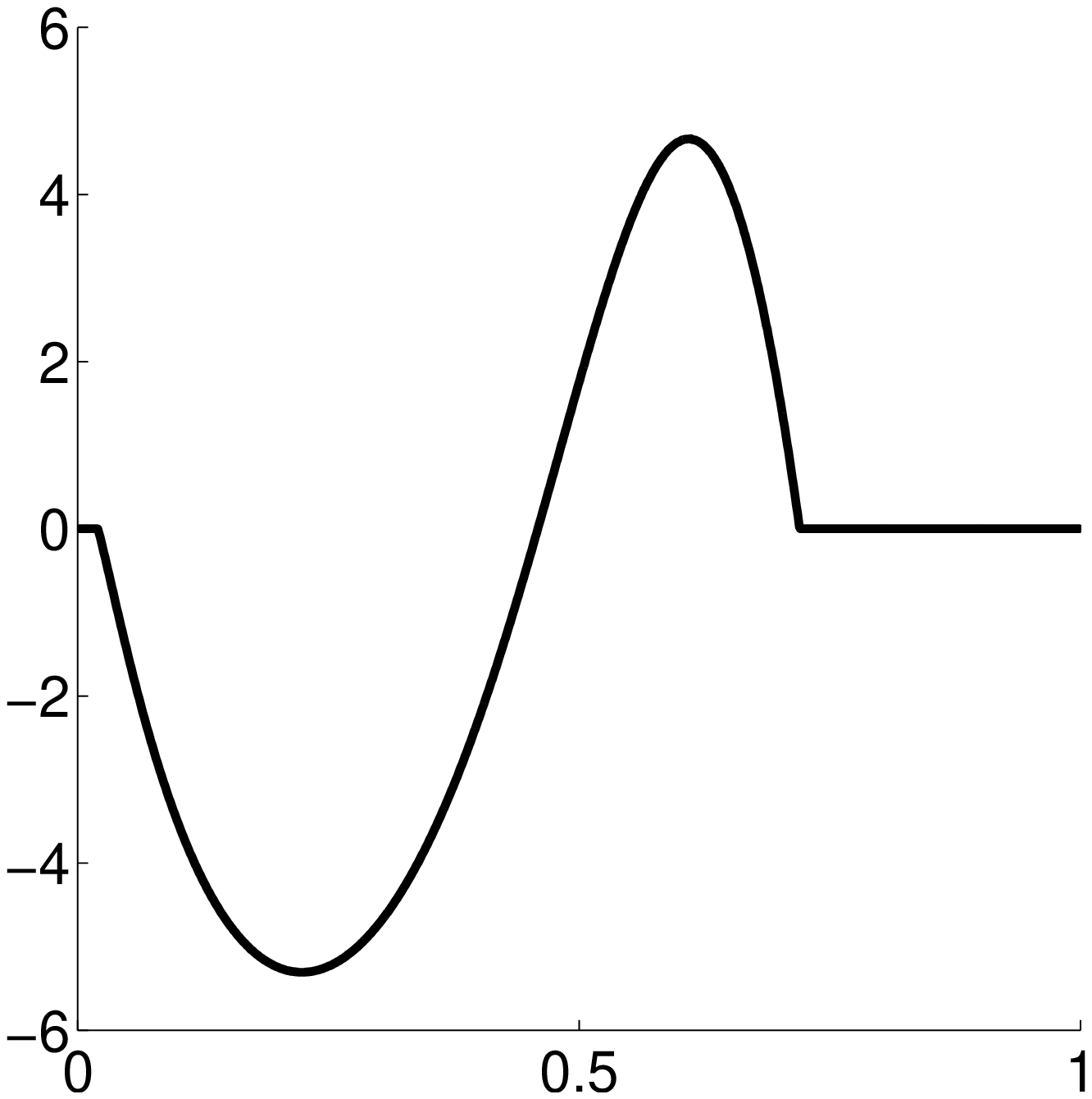}}
\put(15,44){\small $\varphi(|z|)$}
\put(55,44){\small $\sigma(|z|)$}
\put(95,44){\small $q_0^{(2)}(|z|)$}
\put(117,-2){\small $|z|$}
\end{picture}
\caption{\label{fig:profiles}Profile plot of the rotationally symmetric functions.}
\end{figure}

Again, fix $E=-1$. We use exactly the same radially symmetric potentials as in the numerical part of \cite{Music2013} and \cite{deHoop2015}:
\begin{eqnarray}
  q_{\alpha}^{(1)} &=& \alpha \varphi,\label{qalpha1-def}\\
  q_{\alpha}^{(2)} &=& \frac{\Delta\sqrt{\sigma}}{\sqrt{\sigma}}+\alpha\varphi, \label{qalpha2-def}
\end{eqnarray}
where $\alpha\in\R$, $\sigma\in C^2(\Omega)$ with $\sigma\geq c>0$ and $\varphi$ is an approximate test function in $C^2(\Omega)$. See figure \ref{fig:profiles} for the profiles of $\varphi$, $\sigma$ and $q_0^{(2)}$.

\begin{figure}
\begin{picture}(120,100)
\epsfxsize=12cm
\put(0,0){\epsffile{\imagepath 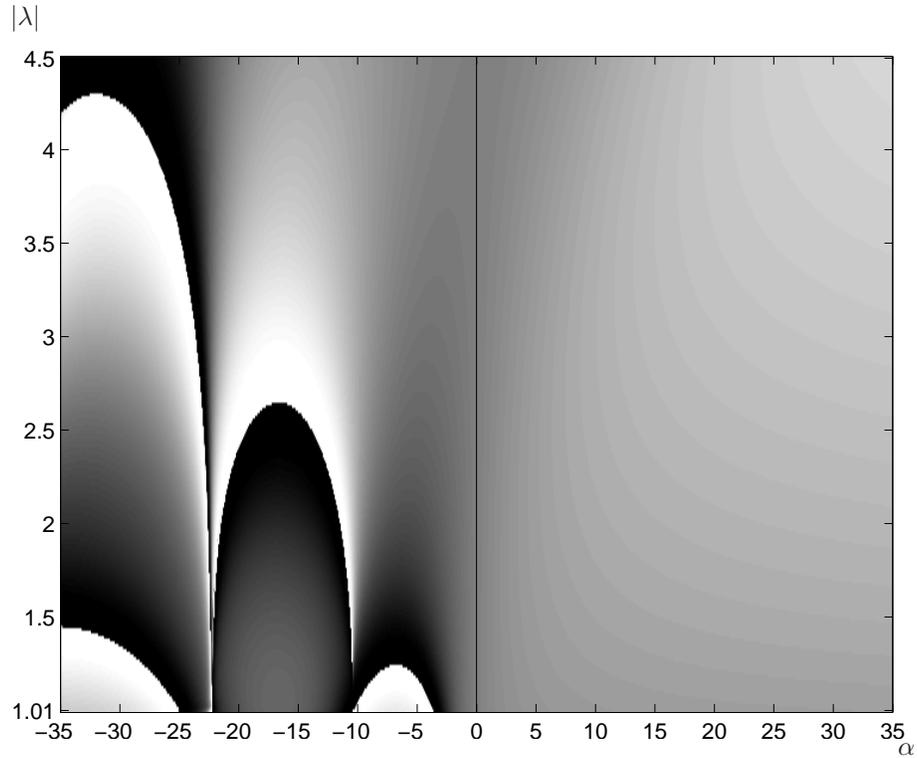}}
\put(118,-2){$\alpha$}
\put(0,95){$\abs{\lambda}$}
\end{picture}
\caption{Scattering transform for the potential $q_\alpha^{(1)} = \alpha\varphi$, energy $E=-1$. x-axis is $\alpha = -35..35$, y-axis is $\lambda = 1.01..4.5$. Compare to figures 3 and 9 in \cite{Music2013} and figure 7 in \cite{deHoop2015}.\label{q1-plane}}
\end{figure}

\begin{figure}
\begin{picture}(120,100)
\epsfxsize=12cm
\put(0,0){\epsffile{\imagepath 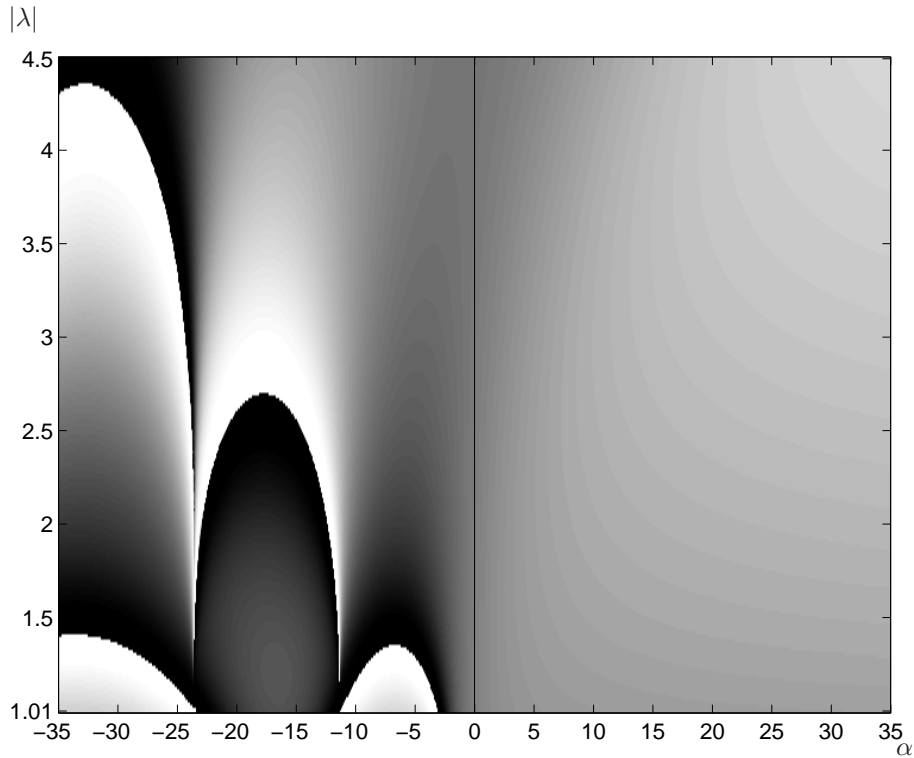}}
\put(118,-2){$\alpha$}
\put(0,95){$\abs{\lambda}$}
\end{picture}
\caption{Scattering transform for the potential $q_\alpha^{(2)} = \Delta\sqrt{\sigma}/\sqrt{\sigma}+\alpha\varphi$, energy $E=-1$. x-axis is $\alpha = -35..35$, y-axis is $\lambda = 1.01..4.5$. Compare to figures 3 and 9 in \cite{Music2013} and figure 7 in \cite{deHoop2015}.\label{q2-plane}}
\end{figure}

We use 250 discrete points of $\lambda$ and 701 discrete points of $\alpha$,
$$
\lambda = 1.01,..,4.5,\quad \alpha = -35,...,35.
$$
For each pair $\{\lambda,\alpha\}$ we compute the CGO solution directly from \eqref{LS-lambda} using the LS solver with $M=8$, leading to $2^M\times 2^M$ sized $z$-grid. Then we compute the radially symmetric and real-valued scattering transform $\vct{t}(\lambda)=\vct{t}(|\lambda|)$ from \eqref{scattering-transform}. In figure \ref{q1-plane} we plot $\vct{t}(|\lambda|)$ for the potential $q_\alpha^{(1)} = \alpha \varphi$; the $x$-axis is the parameter $\alpha$ and the $y$-axis is the modulus $|\lambda|$ of the spectral parameter. In figure \ref{q2-plane} we have the same for  $q_\alpha^{(2)} = \Delta\sqrt{\sigma}/\sqrt{\sigma}+\alpha \varphi$. Black color represents very small negative values, and white very large positive values of $\vct{t}(\lambda)$. The lines where it abruptly changes between these colors are exceptional circles that move as the parameter $\alpha$ changes.

\section{Reconstruction of conductivities at negative energy}\label{sec:conjecture}
In this section we consider potentials of the form 
$$
q = \sigma^{-1/2}\Delta\sigma^{1/2}-E,\qquad E<0,
$$
where $\sigma(z)>0$ is called the \textit{conductivity}, named after the EIT case of $E=0$. Assume that $\sigma|_{\DOm} = s$. Based on numerical evidence, we have
\begin{equation}\label{sigma-recon}
\sigma(z) \approx s\cdot\lim_{\abs{\lambda}\to r^\ast}\re{(\mu(z,\lambda))}^2,
\end{equation}
where $r^\ast$ is an unknown radius. The value $r^\ast\approx 2.5$ gives the best results in our numerical tests. It might be that in reality $r^\ast = 1$ and the computational error in $g_\lambda(z)$ results to this value that is larger than one. Also it could be that in reality there is an integral across the unit circle in \eqref{sigma-recon}, since computing an average over the CGO solutions $\mu(z,\lambda)$ for which $\abs{\lambda}=r^\ast$ also improve the reconstruction. The unit circle $\abs{\lambda}=1$ is special, as seen from the sign function in \eqref{D-bar}. In the positive-energy case the CGO solutions will have a jump when $\lambda$ crosses the unit circle. In the zero-energy case, similar equation to \eqref{sigma-recon} holds, there the spectral parameter is put to zero, see for example \cite{Nachman1996}.

\subsection{Example reconstructions}
\begin{figure}[h]
\begin{picture}(120,125)
\epsfxsize=6cm
\put(1,61){\epsffile{\imagepath 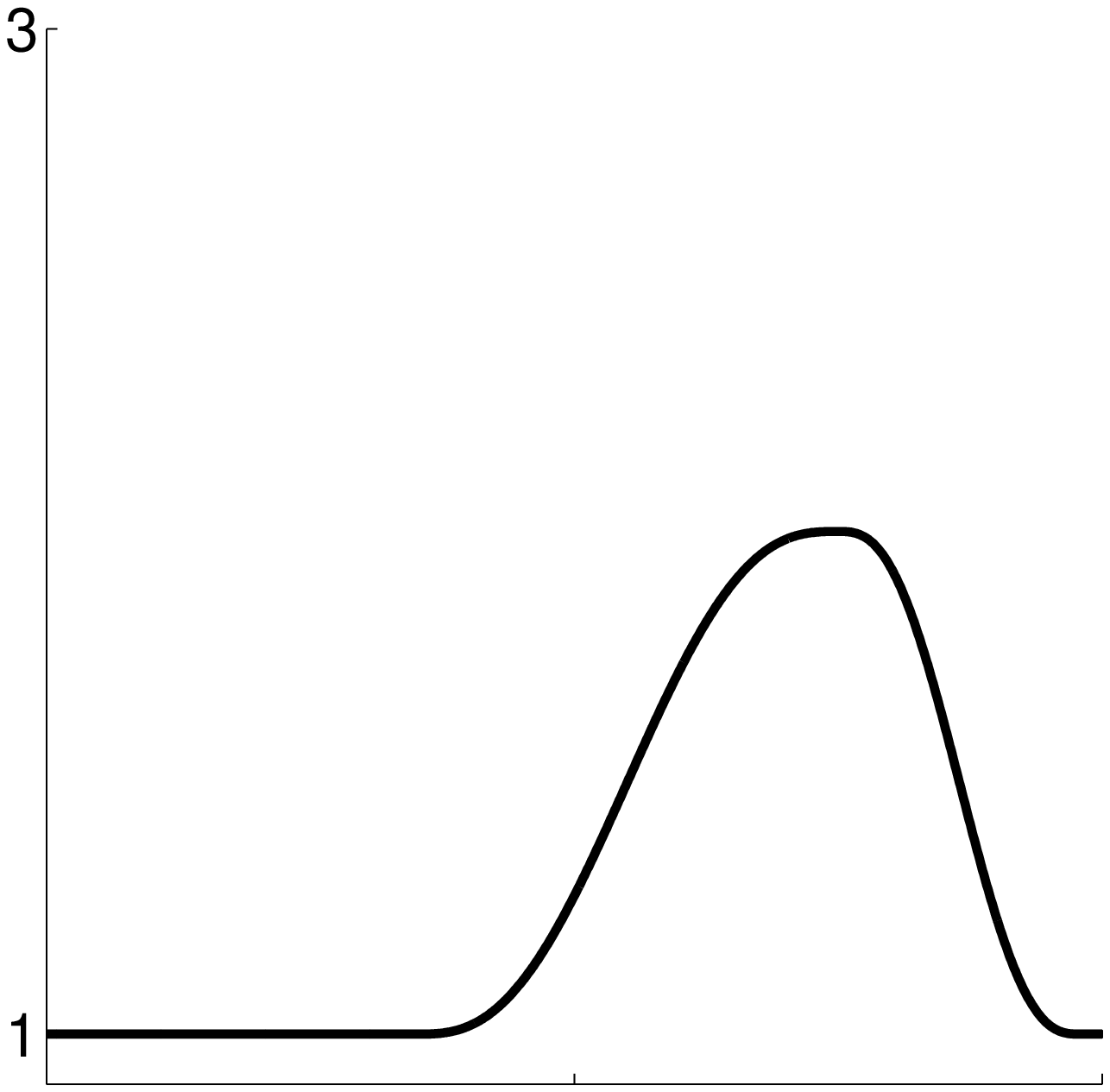}}
\epsfxsize=6cm
\put(61,62){\epsffile{\imagepath 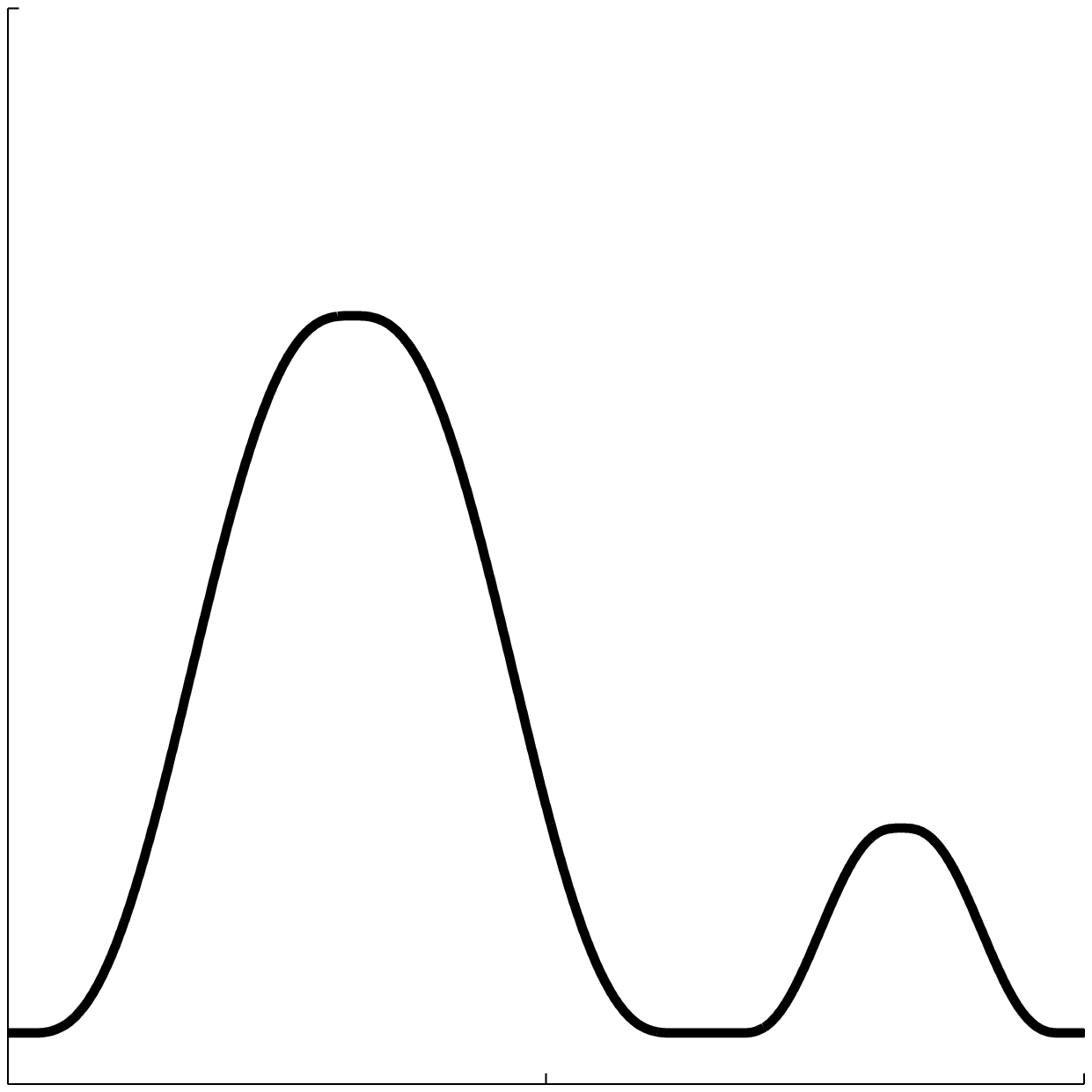}}
\epsfxsize=6cm
\put(0,0){\epsffile{\imagepath 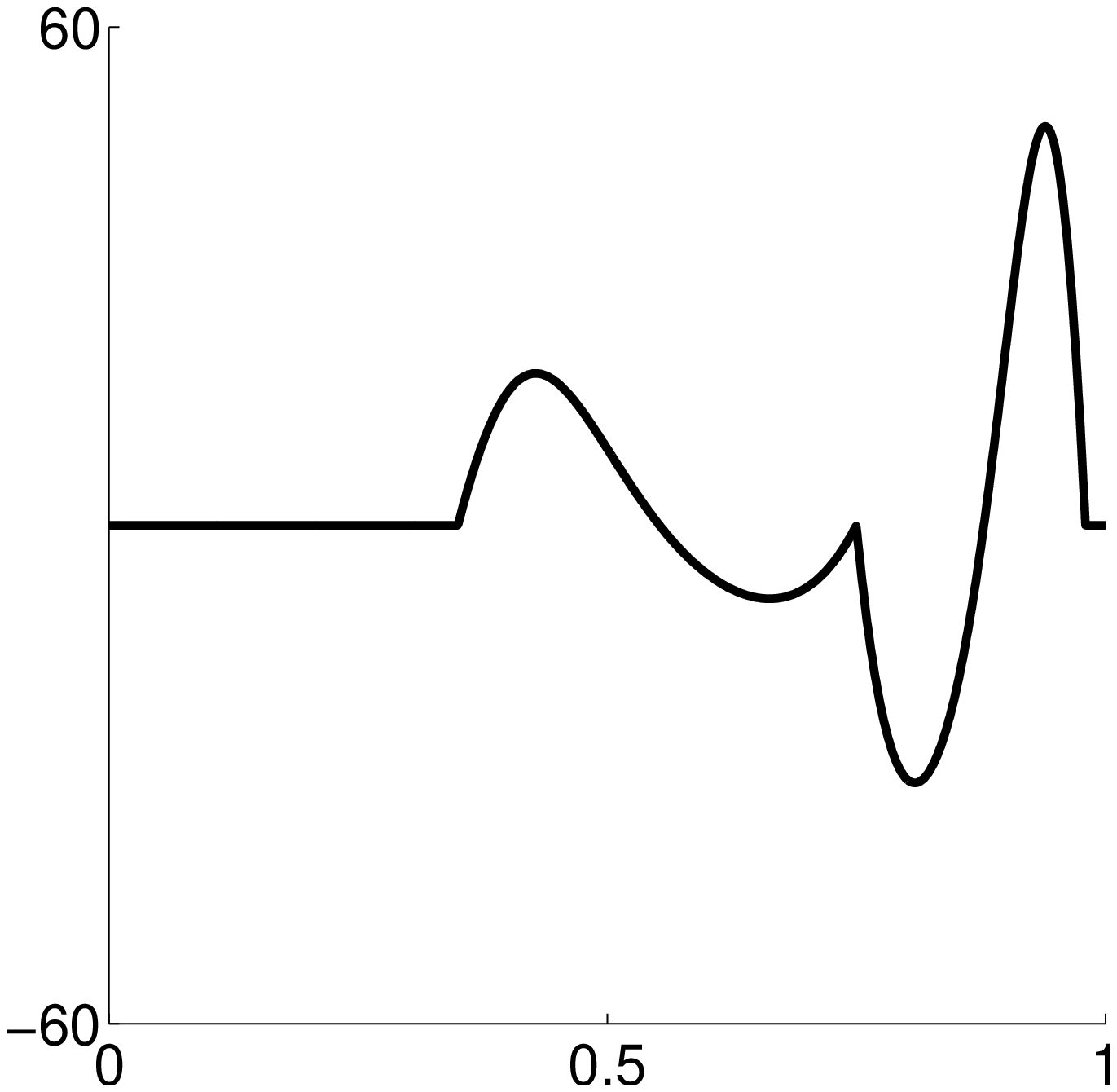}}
\epsfxsize=6cm
\put(60,0){\epsffile{\imagepath 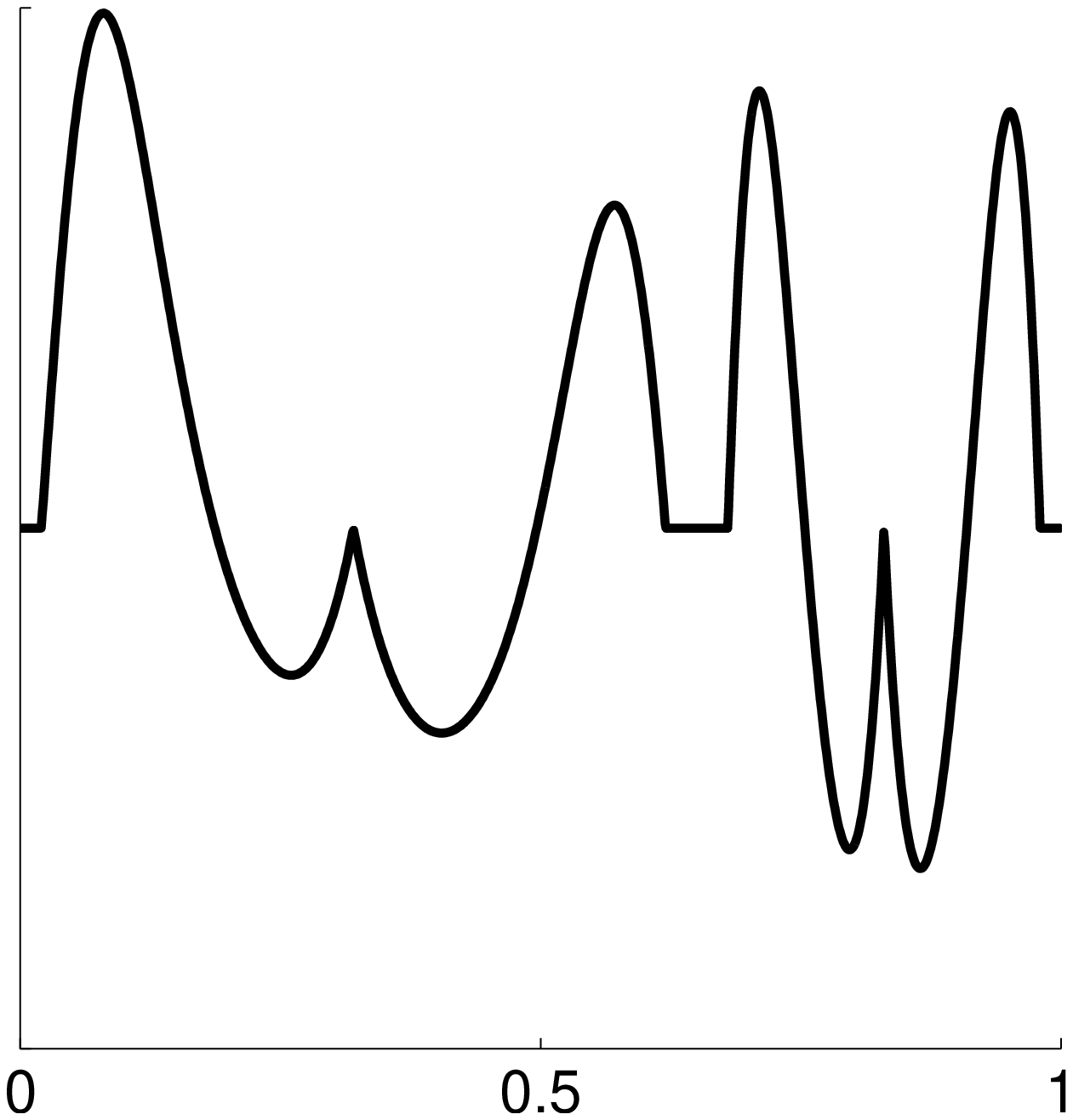}}
\put(30,120){Case 3}
\put(90,120){Case 4}
\put(6,118){$\sigma(\abs{z})$}
\put(6,53){$q_0(\abs{z})$}
\put(119,-3){\small $\abs{z}$}
\end{picture}
\caption{\label{fig:s-profiles}Profile plots of Case 3  and Case 4 radial conductivities $\sigma(z) = \sigma(\abs{z})$ on top. Resulting potentials $q_0 = \sigma^{-1/2}\Delta\sigma^{1/2}$ on the bottom.}
\end{figure}

Again fix $E=-1$. We test the aforementioned method using radial conductivities Case 3 and Case 4 of figure \ref{fig:s-profiles} on top, resulting to radial potentials on the bottom of the figure.

\begin{figure}[h]
\begin{picture}(120,65)
\epsfxsize=6cm
\put(0,1){\epsffile{\imagepath 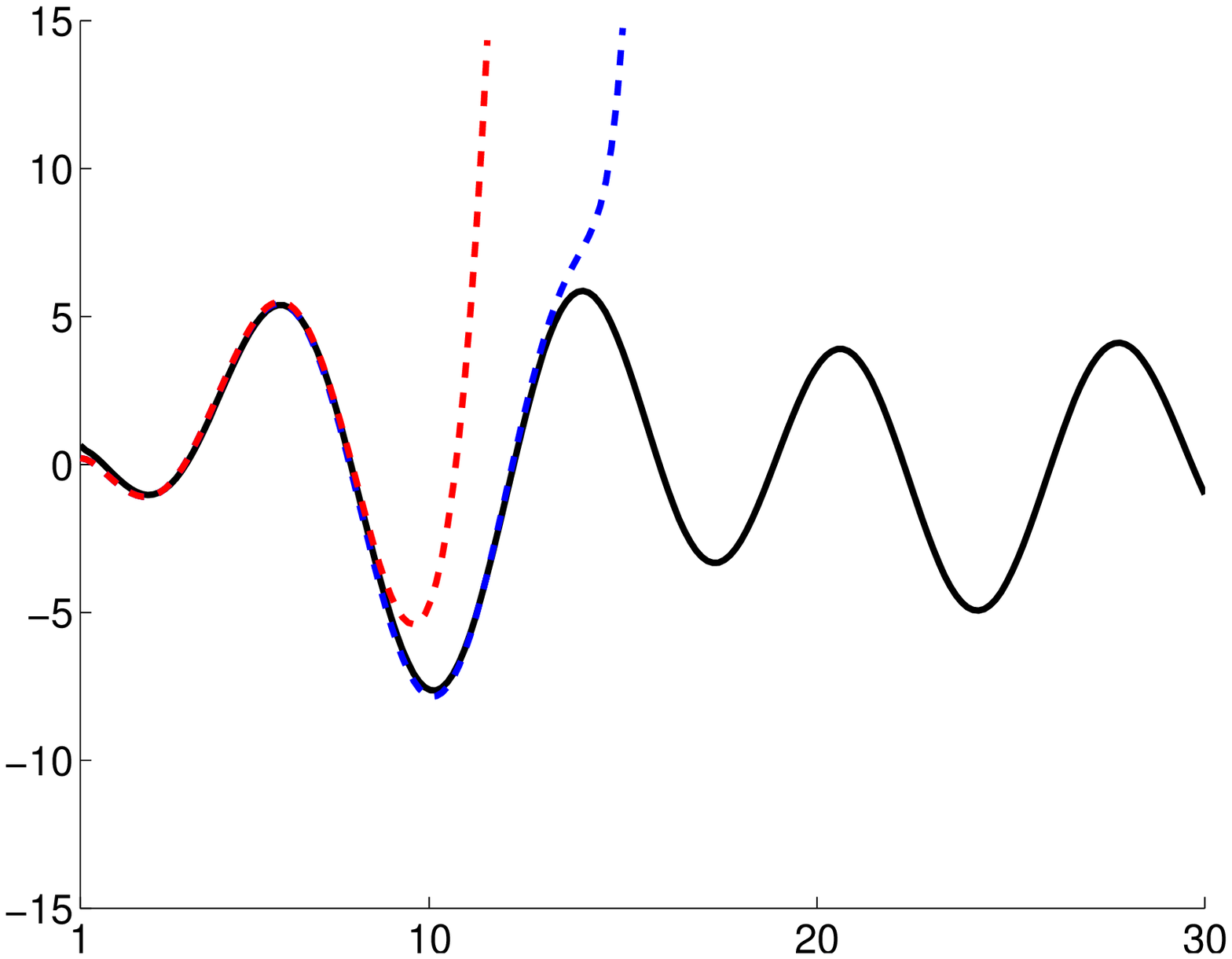}}
\epsfxsize=6cm
\put(60,1){\epsffile{\imagepath 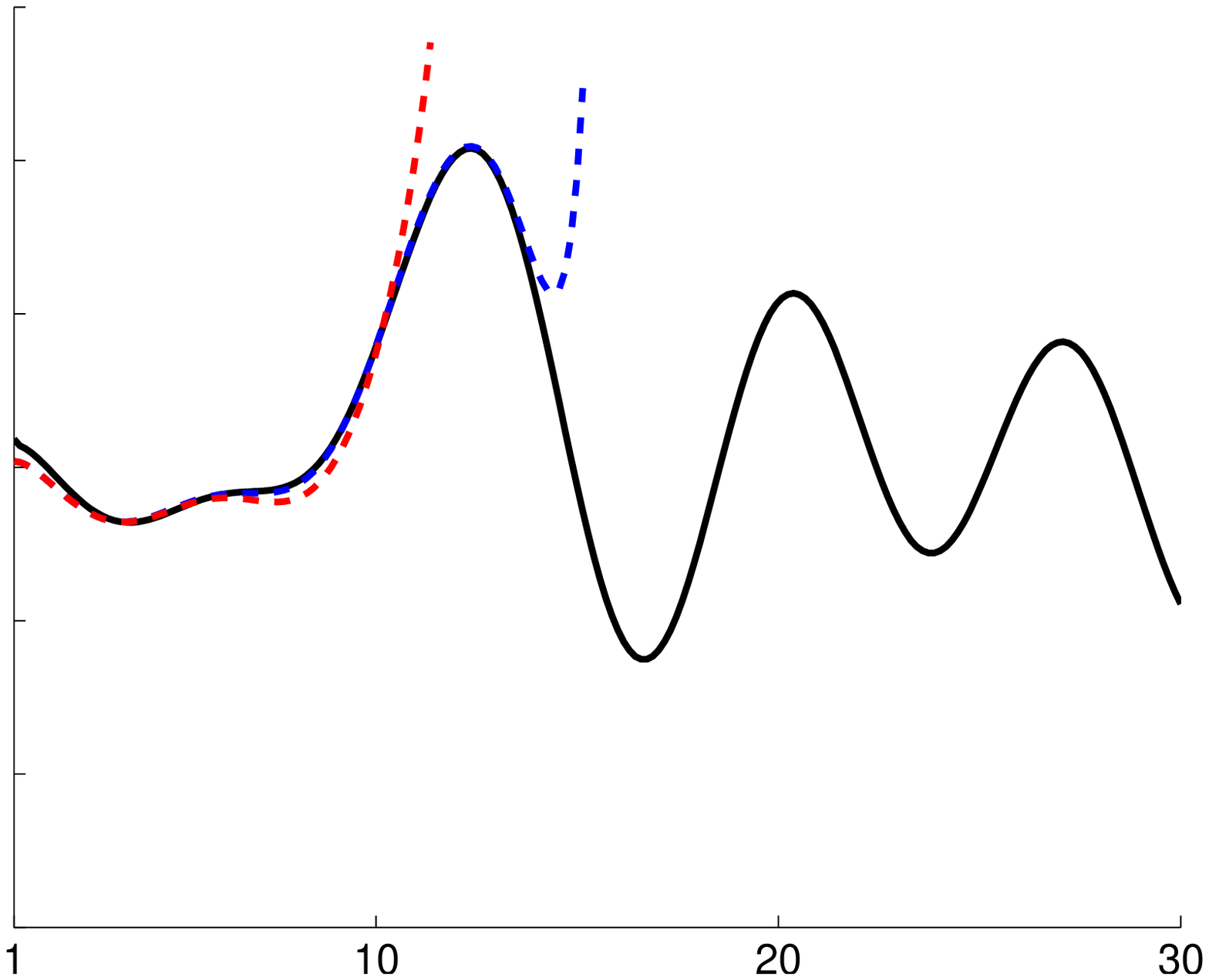}}
\put(28,54){Case 3}
\put(86,54){Case 4}
\put(120,-2){\small $\abs{\lambda}$}
\put(-1,50){\small $\vct{t}(\abs{\lambda})$}
\end{picture}
\caption{\label{fig:pot-scats}The radial scattering transform $\vct{t}(\abs{\lambda})$ for the conductivities of Case 3 and Case 4, at energy $E=-1$, computed in three ways; black solid line indicates computation directly from \eqref{LS-lambda} using the knowledge of $q_0$, blue dashed line indicates computation using the DN-matrices $\mtx{L}_q$, $\mtx{L}_{-E}$ and using equations \eqref{BIE-discrete} and \eqref{scat_matrix}, red dashed line indicates the same but with noisy DN-matrix $\mtx{L}^\epsilon_q$.}
\end{figure}

\begin{figure}[h]
\begin{picture}(120,65)
\epsfxsize=6cm
\put(0,1){\epsffile{\imagepath 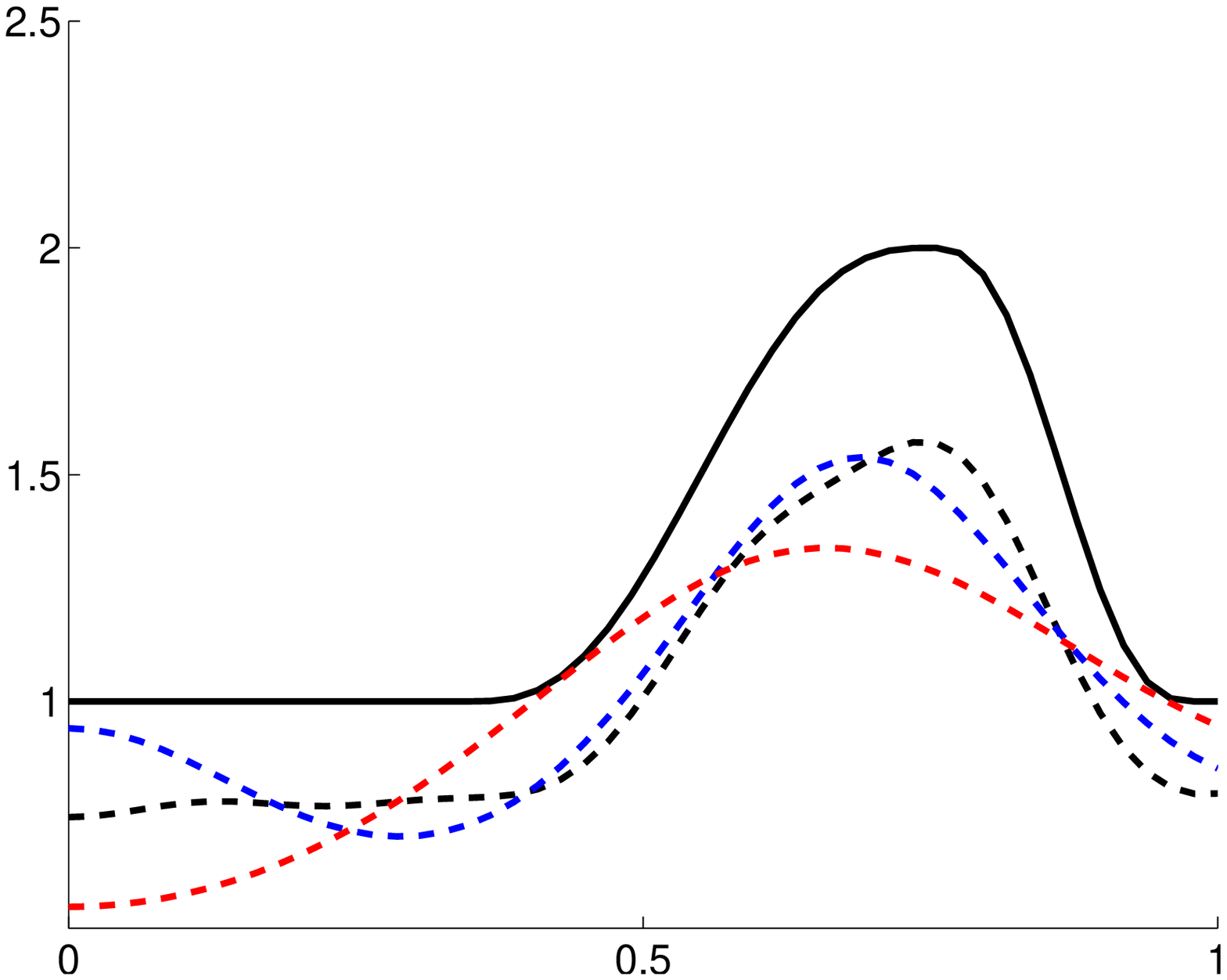}}
\epsfxsize=6cm
\put(60,1){\epsffile{\imagepath 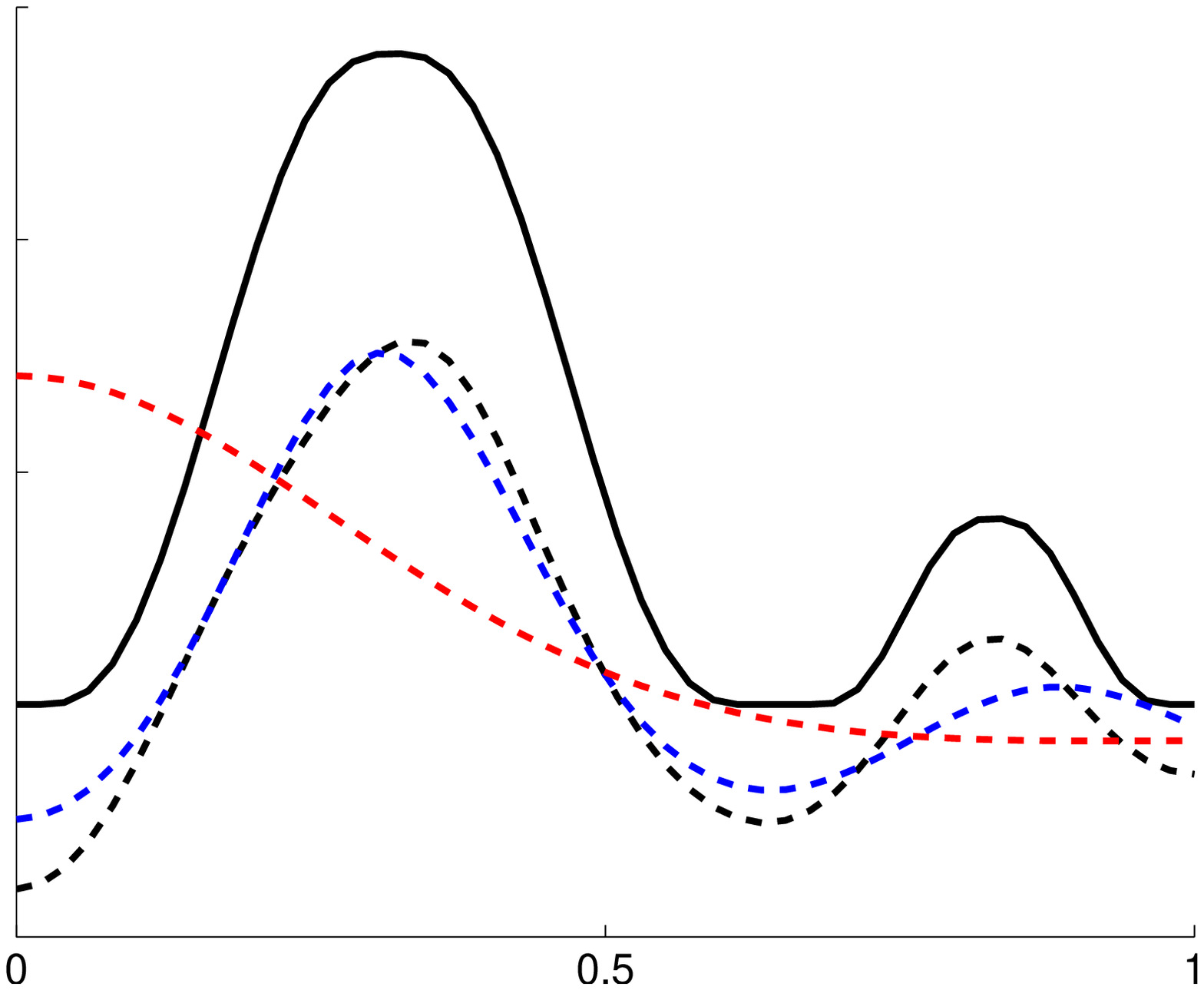}}
\put(26,56){Case 3}
\put(88,56){Case 4}
\put(-1,52){\small $\sigma(\abs{z})$}
\put(120,-2){\small $\abs{z}$}
\end{picture}
\caption{\label{fig:s-recon}Reconstructions of test conductivities Case 3 and Case 4, at energy $E=-1$. Black solid line: the original conductivity. Black dashed line: reconstruction using the knowledge of the potential via the scattering transform \eqref{scattering-transform}. Blue dashed line: reconstruction using the new computational method and \eqref{sigma-recon} without noise. Red dashed line: the same but with added noise.}
\end{figure}

The radial scattering transforms are pictured in figure \ref{fig:pot-scats}. The reconstructions of figure \ref{fig:s-recon} are computed by equation \eqref{sigma-recon} with $r^\ast = 2.5$ and by computing an average of $\mu(z,\lambda)$ for which $\abs{\lambda}=r^\ast$. 


\section{Application to Diffuse Optical Tomography}\label{sec:DOT}
\begin{figure}
\begin{picture}(120,55)
\epsfxsize=4cm
\put(0,0){\epsffile{\imagepath 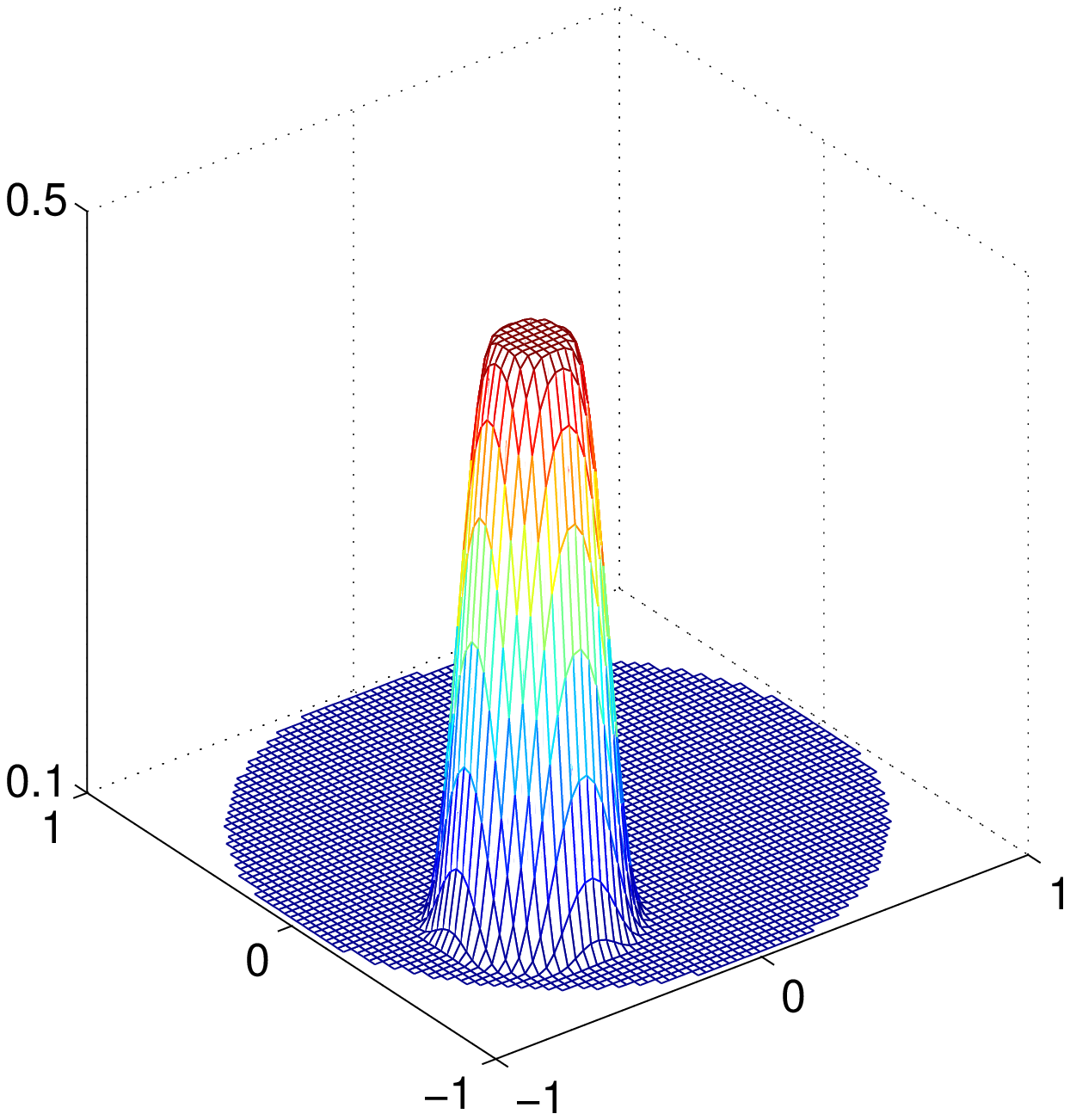}}
\epsfxsize=4cm
\put(40,0){\epsffile{\imagepath 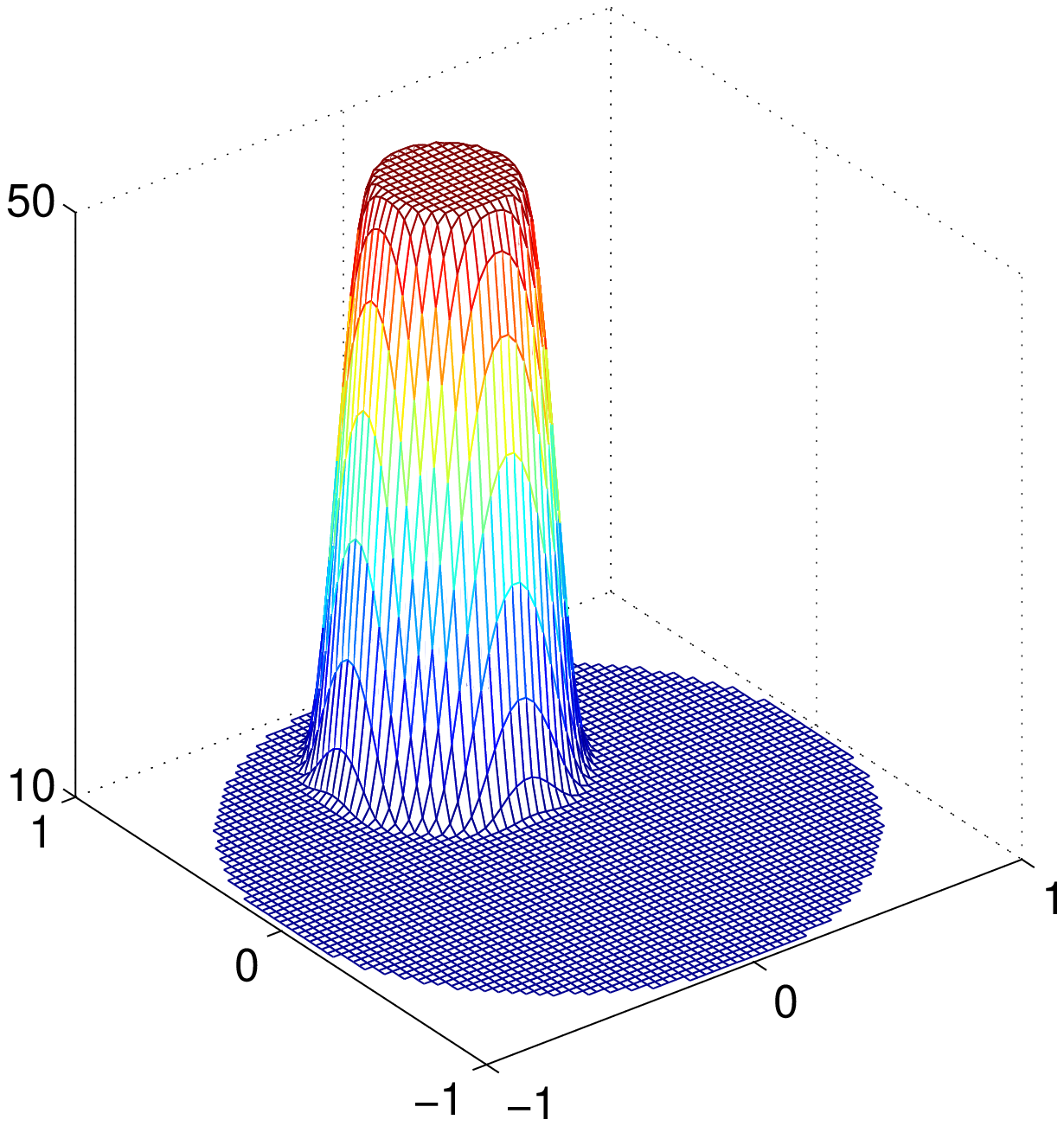}}
\epsfxsize=4cm
\put(80,0){\epsffile{\imagepath 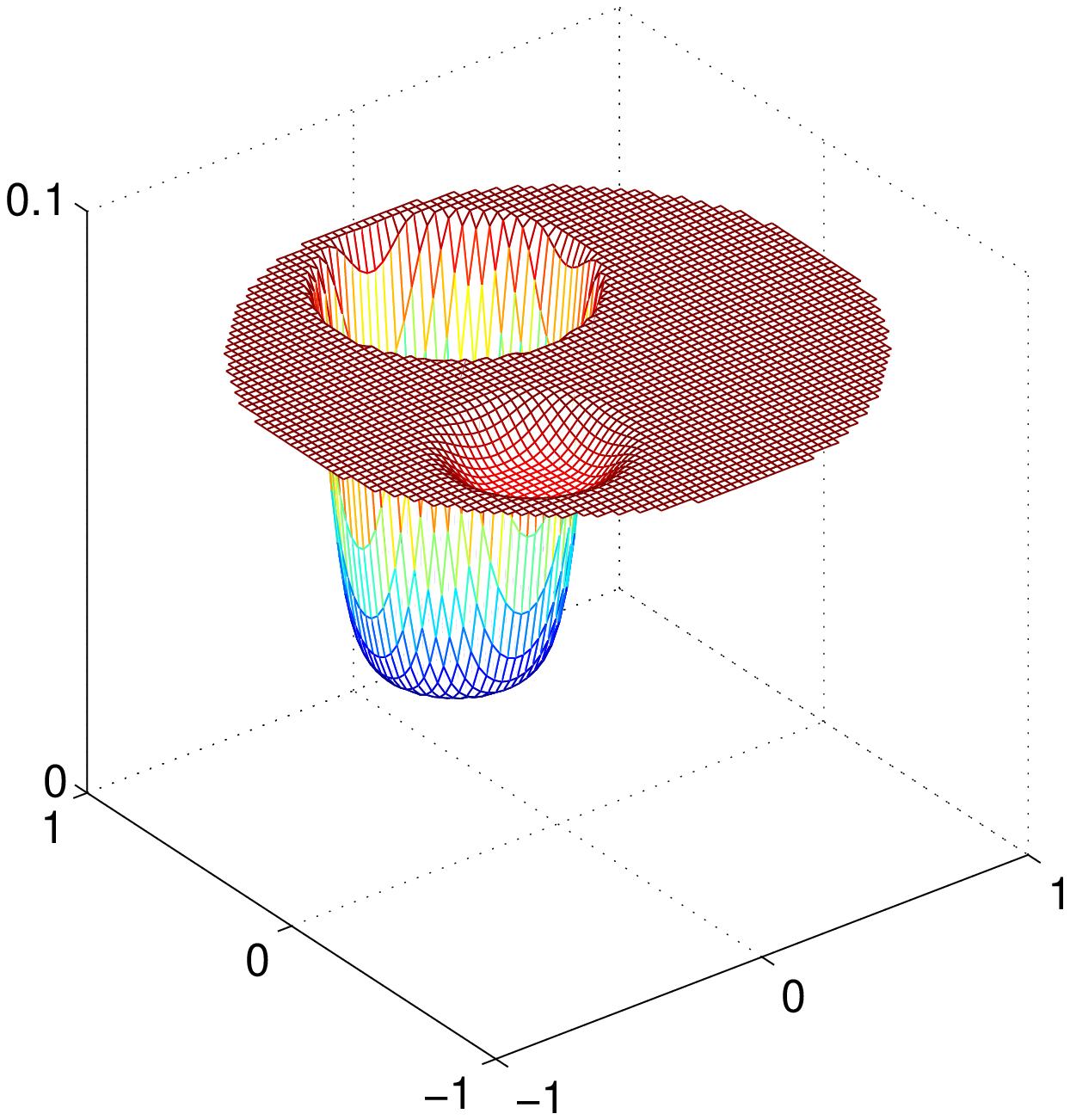}}
\put(10,45){$\mu_a(z)$, 1/cm}
\put(50,45){$\mu_s(z)$, 1/cm}
\put(97,45){$D(z)$}
\end{picture}
\caption{\label{fig:DOTcoeff} The DOT setting: absorption coefficient on the left, scattering coefficient in the middle and the resulting diffusion coefficient on the right.}
\end{figure}

\begin{figure}
\begin{picture}(120,67)
\epsfxsize=5.5cm
\put(0,0){\epsffile{\imagepath 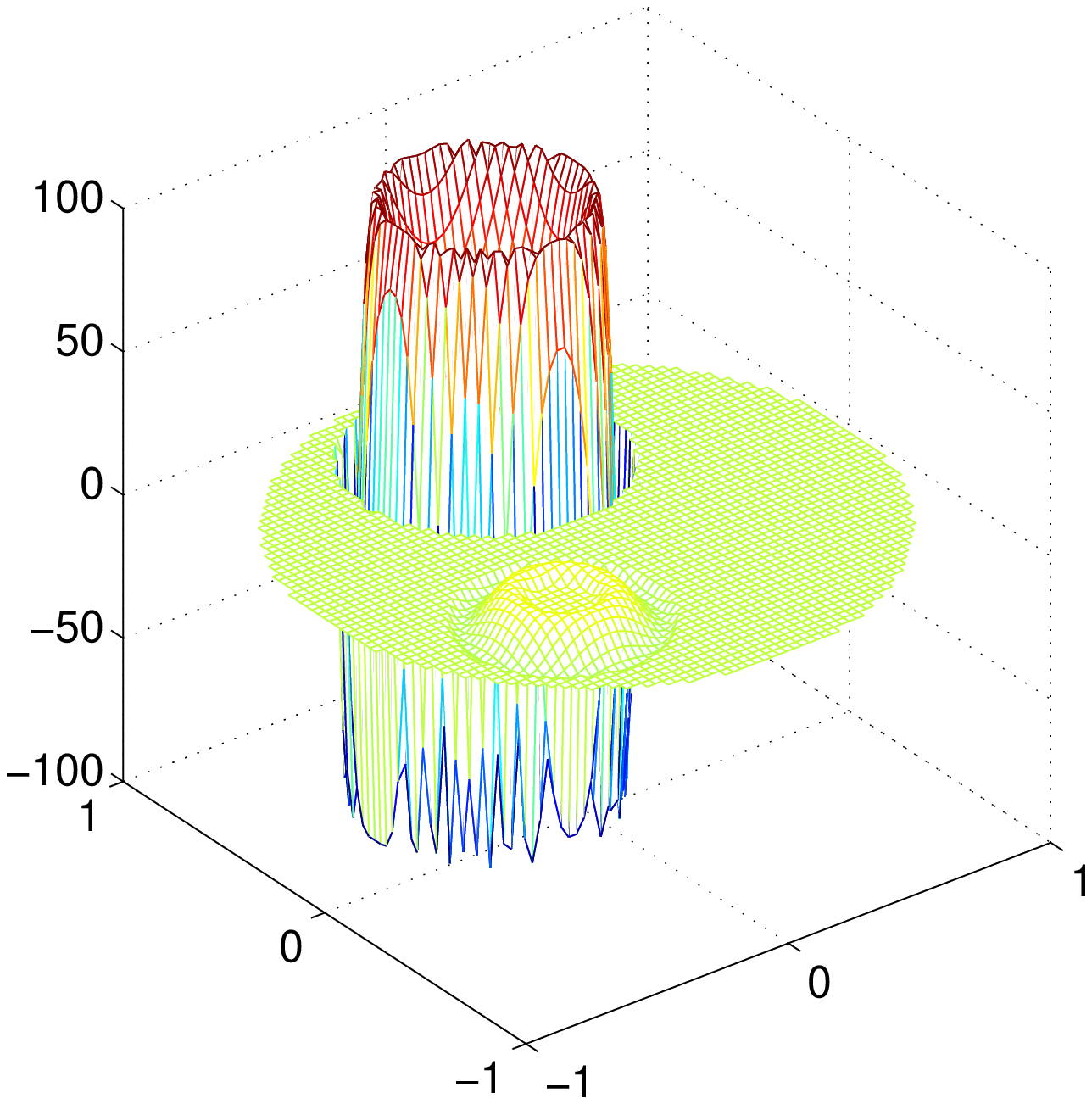}}
\epsfxsize=5.5cm
\put(60,0){\epsffile{\imagepath 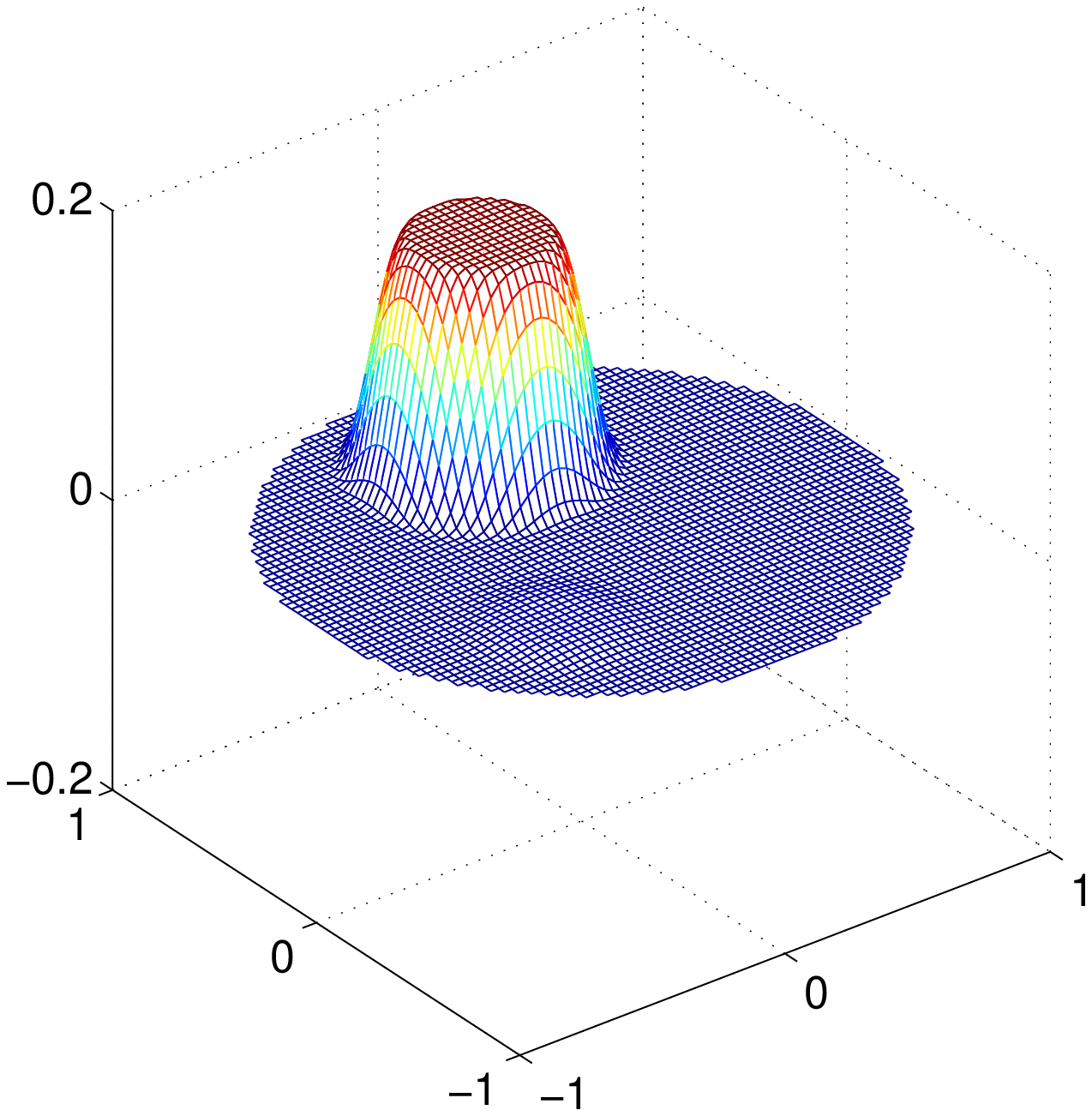}}
\put(23,59){$\re(q_0(z))$}
\put(71,59){$\im(q_0(z))$, $100\MHz$}
\end{picture}
\caption{\label{fig:DOTq0}The potential $q_0 = D^{-1/2}\Delta D^{1/2}+\frac{1}{D}(\mu_a+\frac{i\omega}{c})$.}
\end{figure}

Recall the DOT problem of section \ref{DOT}. In addition to the physical parameters presented there we also have the scattering anisotropy parameter $g$, the scattering coefficient $\mu_s$ and the relations
$$
\mu_s^\prime = (1-g)\mu_s,\quad D = \frac{1}{3(\mu_a+\mu_s^\prime)}.
$$
\begin{itemize}
\item Absorption coefficients $\mu_a$ typically vary between 0.1 and 0.5 1/cm. This scaling leads to the unit disc corresponding to the radius of 1 cm.
\item Scattering coefficients $\mu_s$ typically vary between 10 and 50 1/cm.
\item For skull and grey matter we have $g=0.6$.
\item We assume angular modulation of $\omega=100 \MHz$ and $\omega=0$.
\end{itemize}
These values motivate our examples. In figure \ref{fig:DOTcoeff} we have the absorption coefficient on the left, the scattering coefficient in the middle and resulting diffusion coefficient on the right. With $\omega=100 \MHz$ this leads to the potential $q_0 = D^{-1/2}\Delta D^{1/2}+\frac{1}{D}(\mu_a+\frac{i\omega}{c})$ of figure \ref{fig:DOTq0} and $E =  -1.23 - 0.041i$. With $\omega=0$ the potential $q_0$ is real-valued and $E=-1.23$.

The method of reconstructing the conductivity at negative energy described in section \ref{sec:conjecture} does not cover neither of these cases because of the extra term $\mu_a/D$ in the potential. Regardless, we proceed to test the method to both cases of $\omega$, where the value $\omega=100 \MHz$ requires an additional approximation since in our method we assume real-valued $q$ and thus negative energy $E<0$. 

\begin{figure}[h]
\begin{picture}(120,128)
\epsfxsize=5.5cm
\put(0,65){\epsffile{\imagepath 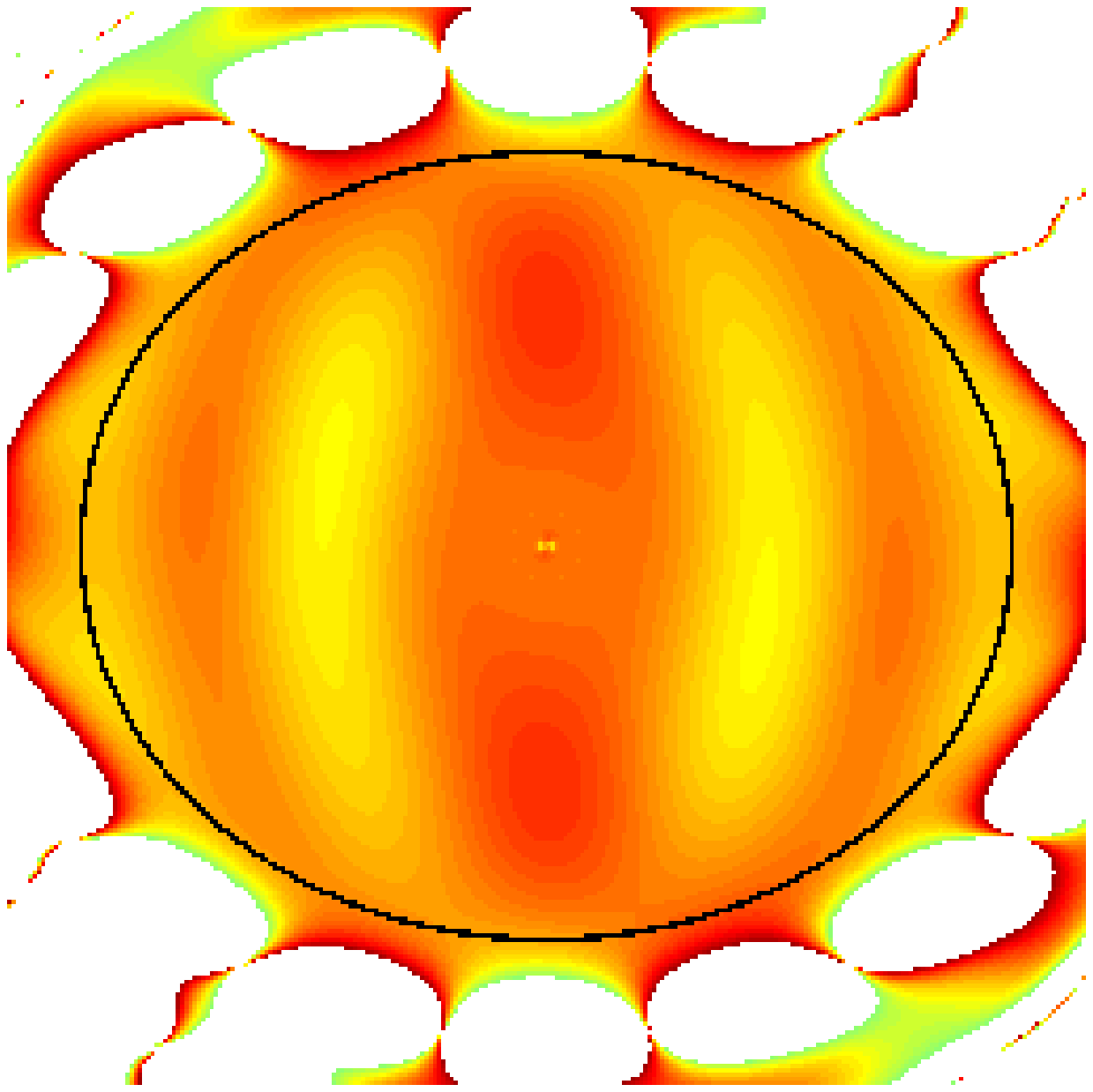}}
\epsfxsize=5.5cm
\put(60,65){\epsffile{\imagepath 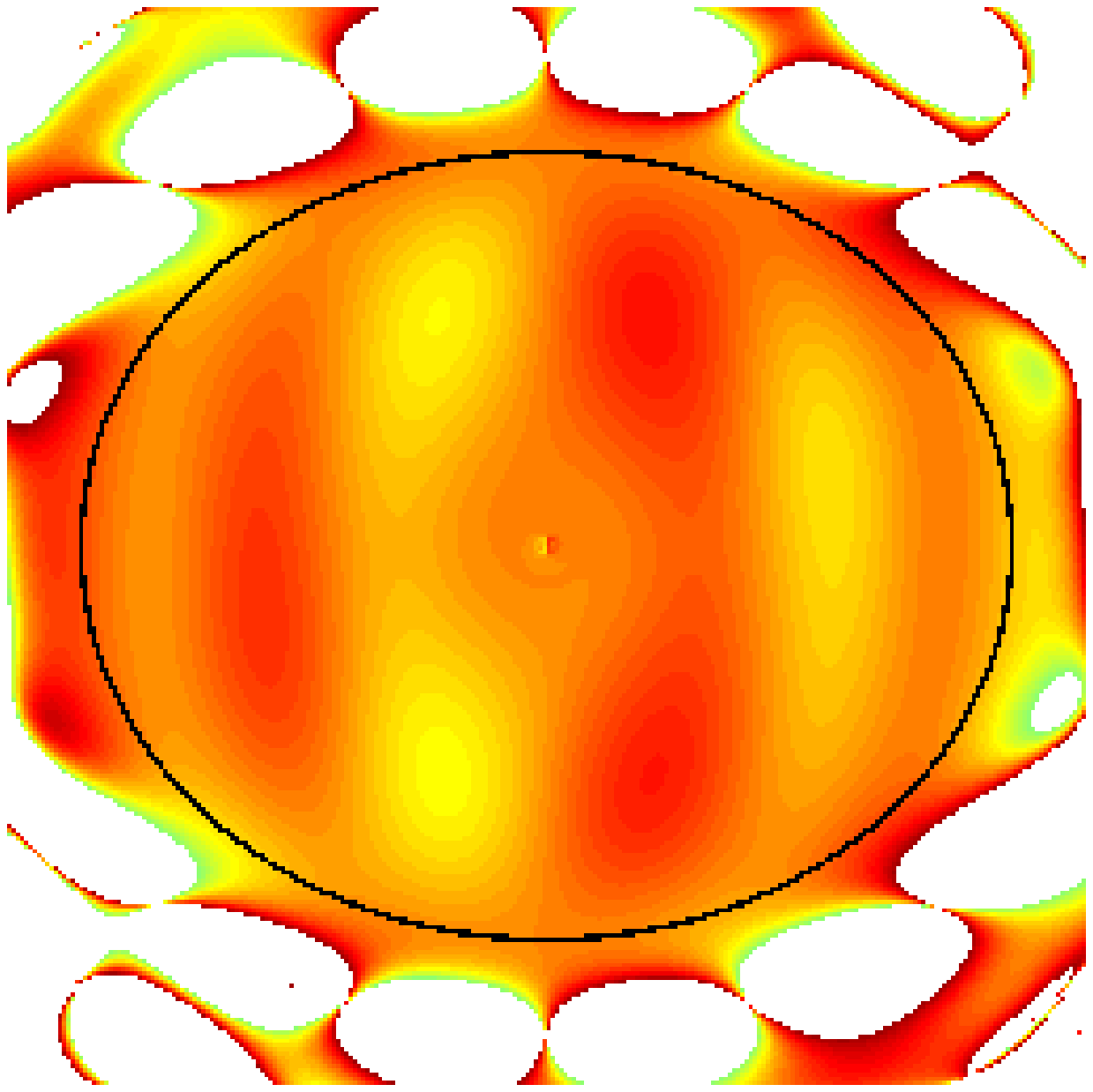}}
\epsfxsize=5.5cm
\put(0,0){\epsffile{\imagepath 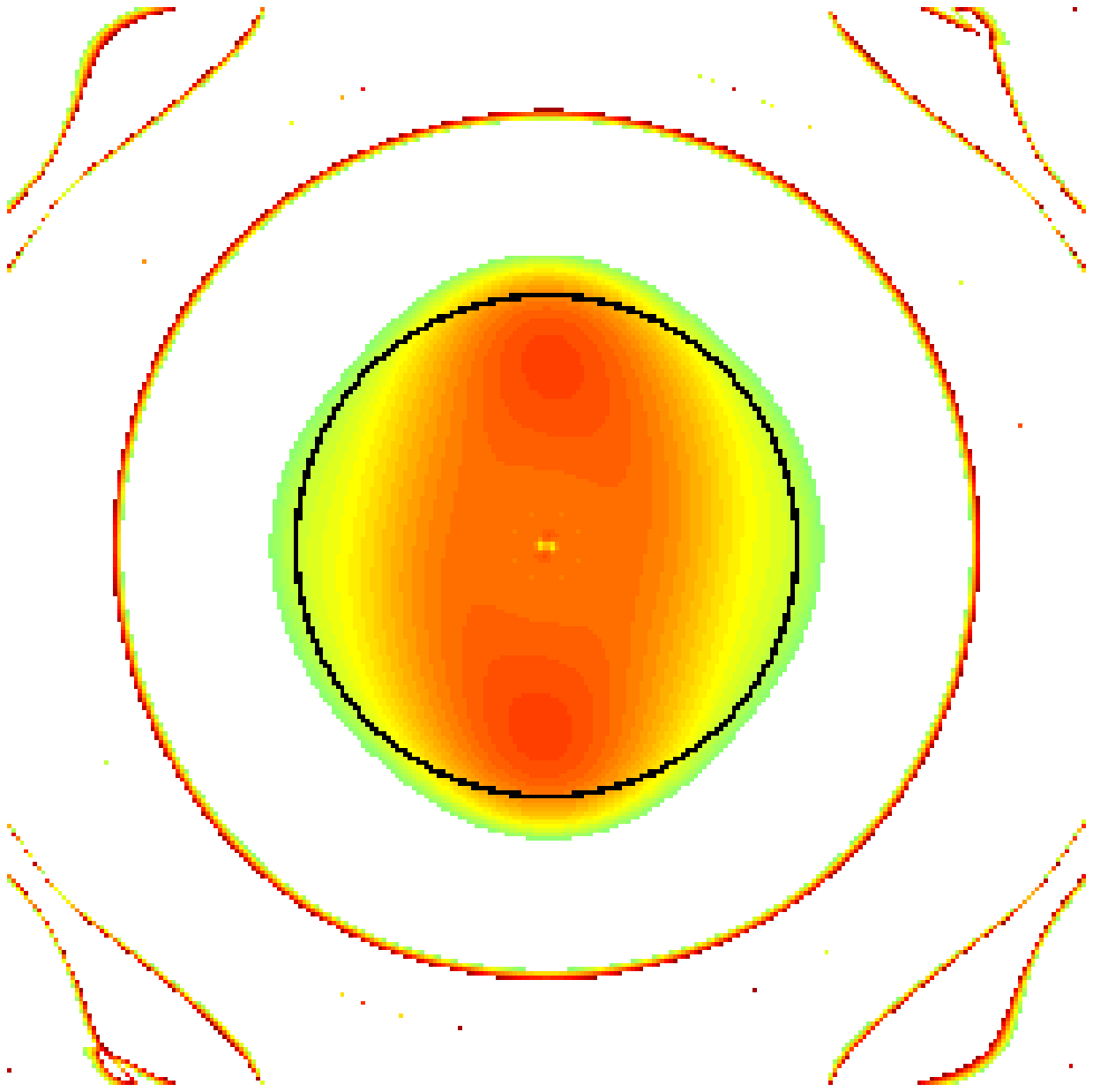}}
\epsfxsize=5.5cm
\put(60,0){\epsffile{\imagepath 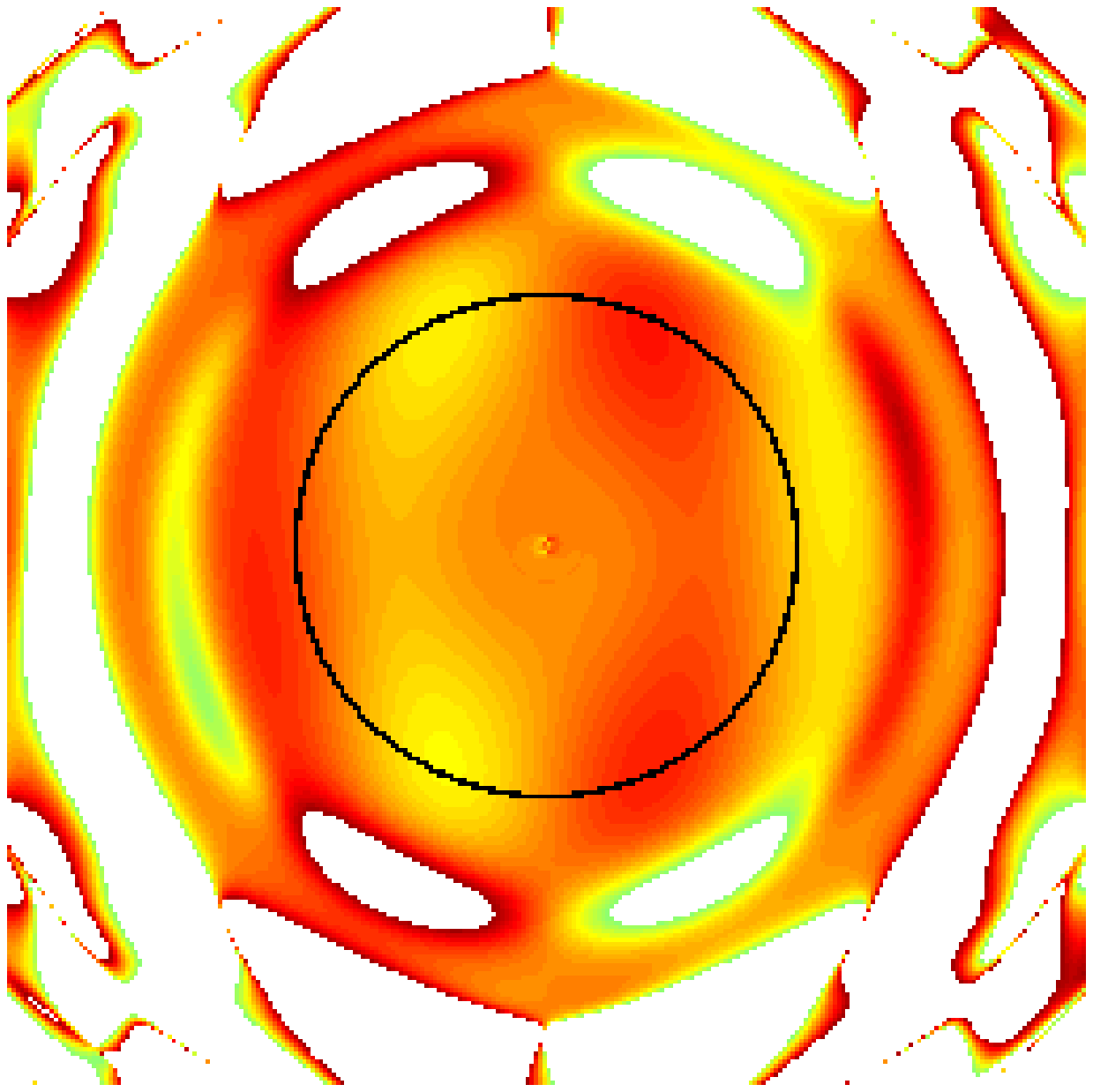}}
\put(21,124){$\re(\vct{t}(\lambda))$}
\put(85,124){$\im(\vct{t}(\lambda))$}
\put(48,124){Non-noisy, $\mathrm{L}_q$}
\put(48,58){Noisy, $\mathrm{L}^{\epsilon}_q$}
\end{picture}
\caption{\label{fig:DOTscat}The scattering transform $\vct{t}(\lambda)$ of the DOT potential of figure \ref{fig:DOTq0} with $\omega=100\MHz$. Real part on the left, imaginary part on the right, in a $\lambda$ -grid $[-15,15]\times [-15,15]i$. On the top row: the non-noisy DN-matrix $\mathrm{L}_q$ was used. On the bottom row: the noisy DN-matrix $\mathrm{L}^{\epsilon}_q$ was used. In the white areas the computation breaks down. The black line indicates the ellipse used for the truncation $\vct{t}_R(\lambda)$.}
\end{figure}

In figure \ref{fig:DOTscat} we have the scattering transform of the potential of figure \ref{fig:DOTq0}, $\omega = 100\MHz$, real part on the left, imaginary part on the right, computed from the non-noisy DN-matrix on top and from noisy DN-matrix on the bottom. In white areas the computation is failing due to numerical errors caused by large values of $\lambda$ and/or noise. Recall the ellipse \eqref{ellipse} used to truncate the scattering transform. The black line is the ellipse used for each case, we used
\begin{itemize}
\item for $\omega = 100\MHz$ without noise $a=11$, $b=13$, $\phi=\pi/2$,
\item for $\omega = 100\MHz$ with noise $a=7$, $b=7$, $\phi=0$,
\item for $\omega = 0\MHz$ without noise $a=11$, $b=13$, $\phi=\pi/2$,
\item for $\omega = 0\MHz$ with noise $a=8$, $b=8$, $\phi=0$.
\end{itemize}
The scattering transform in the case $\omega=0$ looks similar to figure \ref{fig:DOTscat}.

\begin{figure}
\begin{picture}(120,160)
\epsfxsize=5cm
\put(10,100){\epsffile{\imagepath 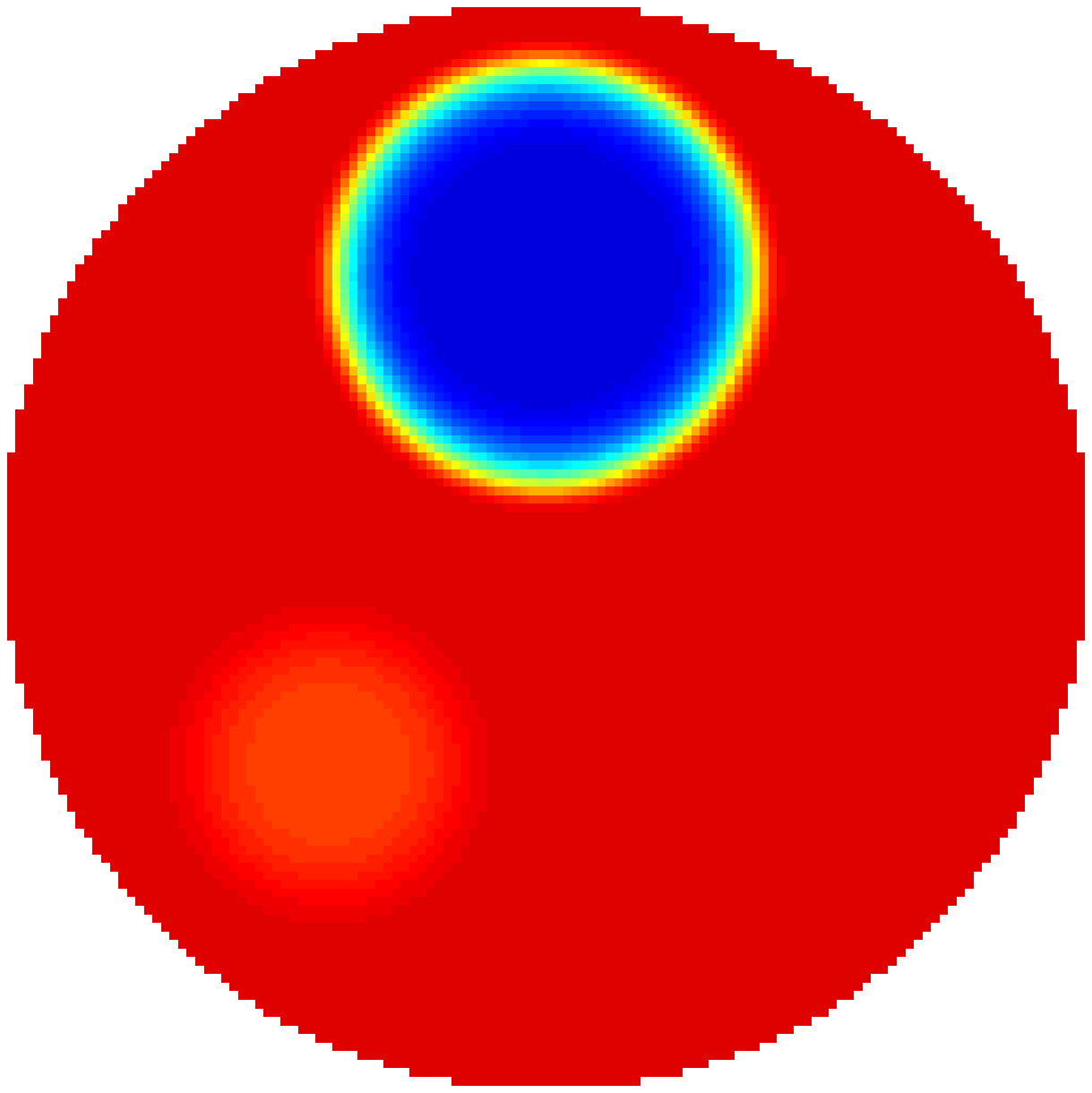}}
\epsfxsize=5cm
\put(60,100){\epsffile{\imagepath DOTv1_Dorig.eps}}
\epsfxsize=5cm
\put(10,50){\epsffile{\imagepath 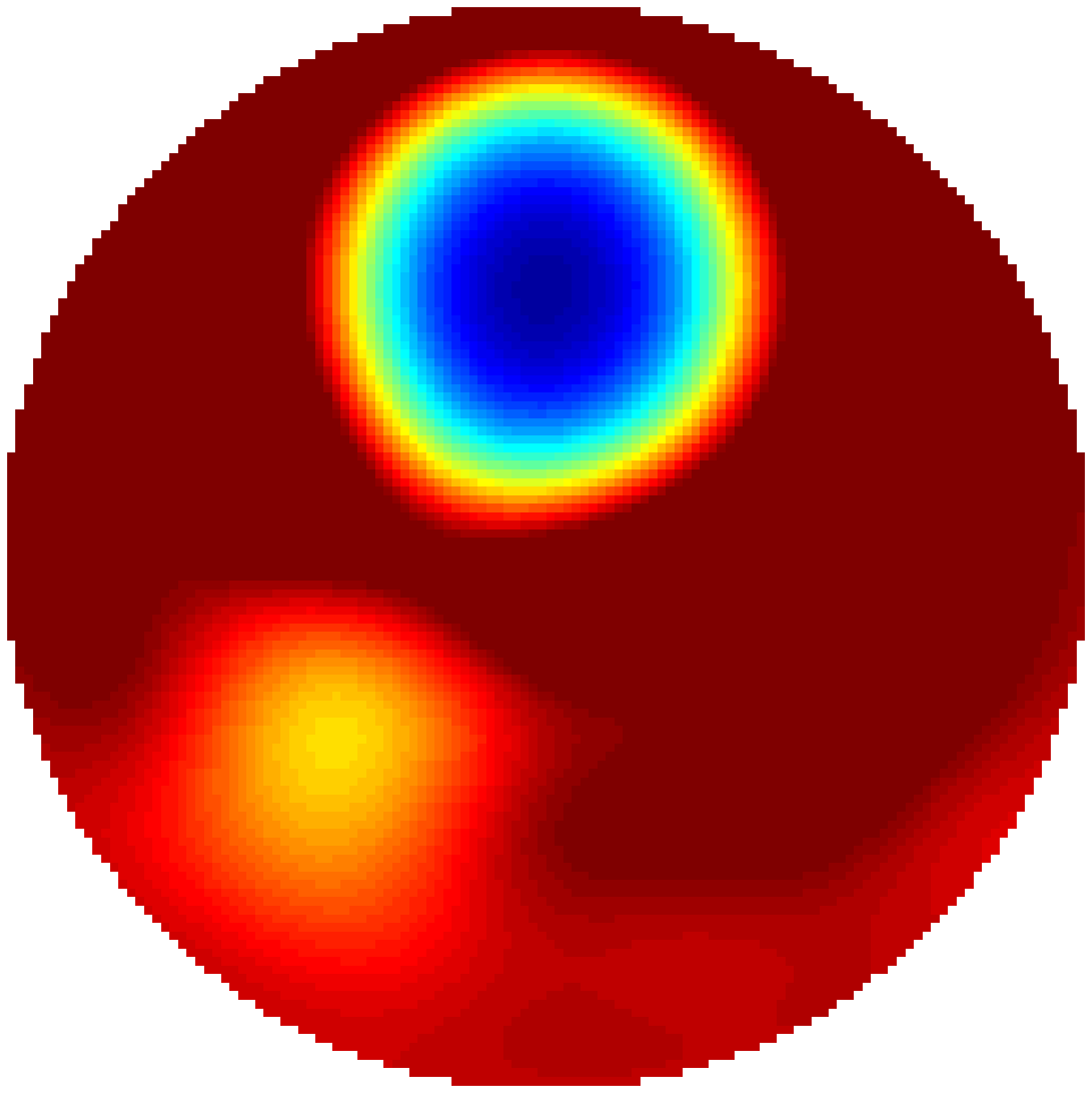}}
\epsfxsize=5cm
\put(60,50){\epsffile{\imagepath 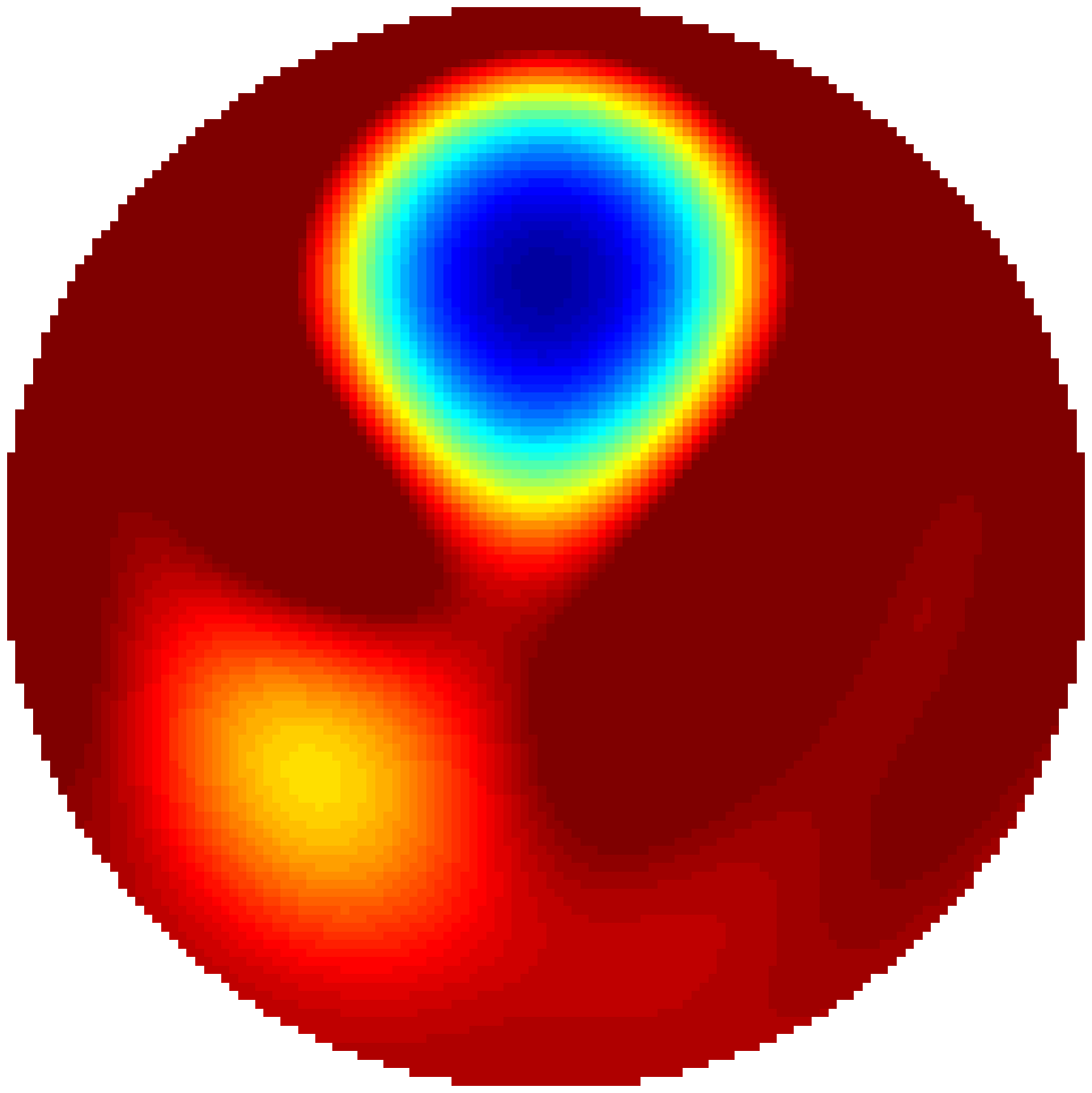}}
\epsfxsize=5cm
\put(10,0){\epsffile{\imagepath 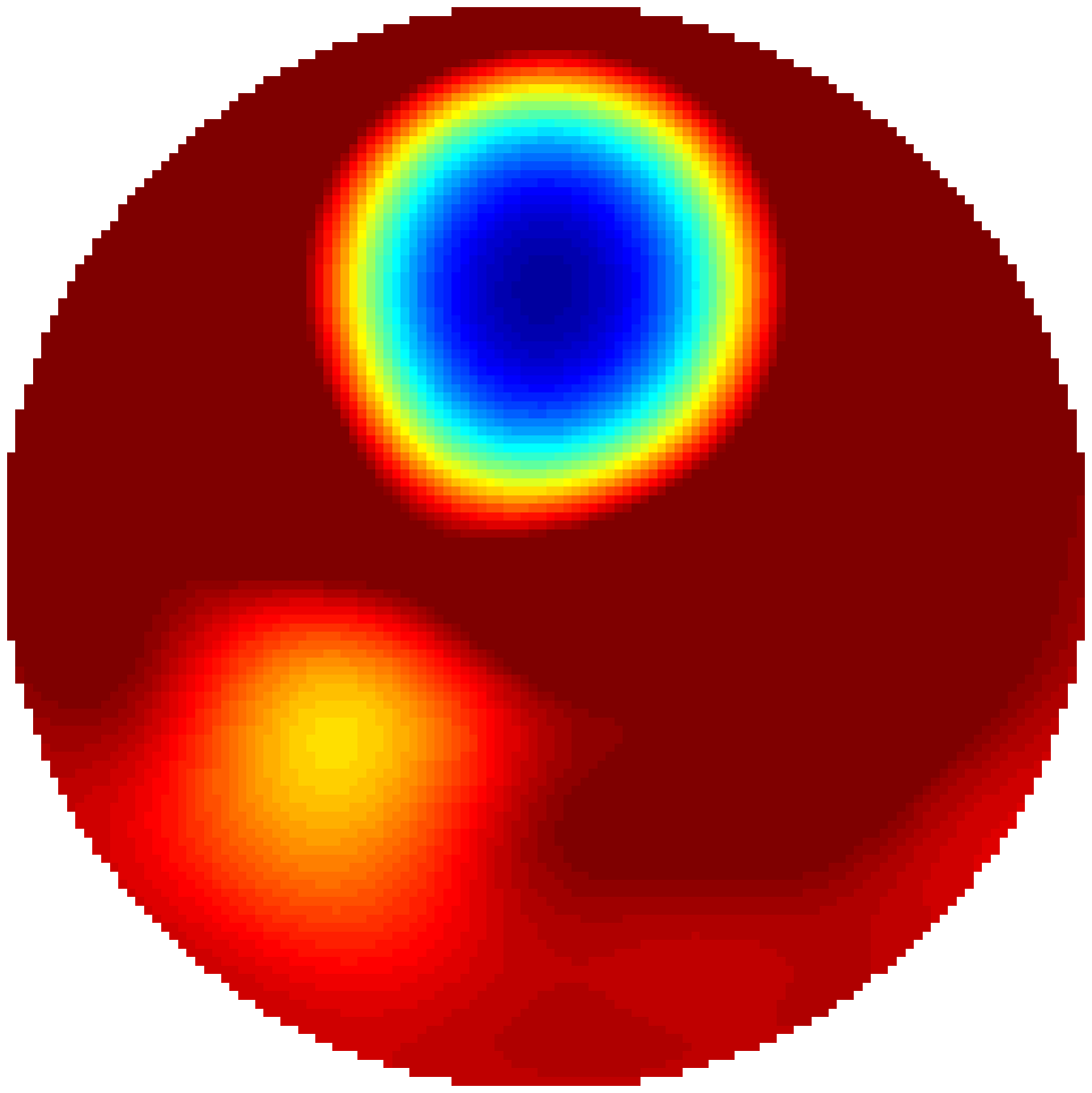}}
\epsfxsize=5cm
\put(60,0){\epsffile{\imagepath 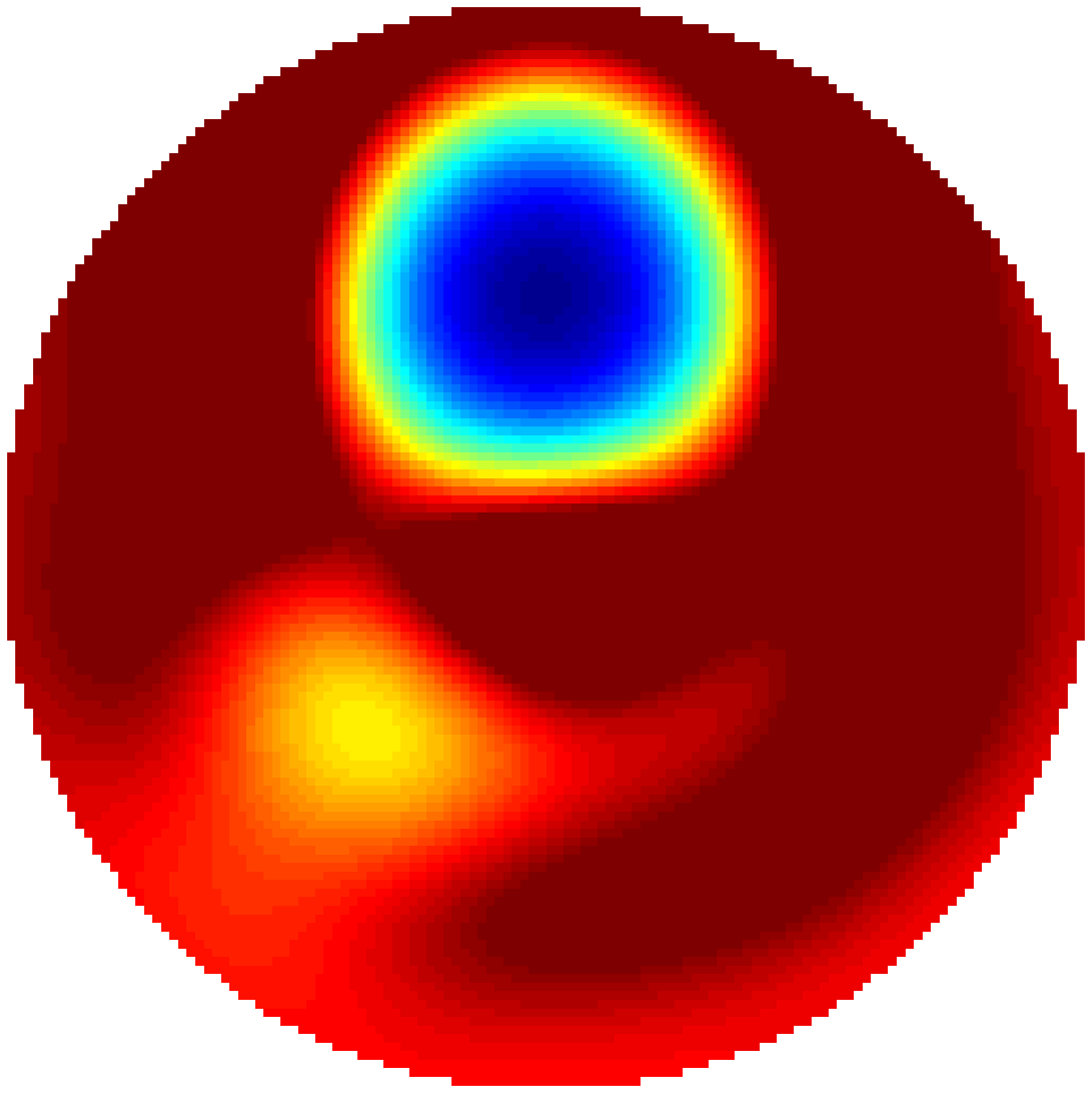}}
\put(21,153){$\omega = 100\MHz$}
\put(74,153){$\omega = 0\MHz$}
\put(53,148){\small Original}
\put(56,144){\small $D(z)$}
\put(50,98){\small Reconstruction}
\put(58,94){\small $\mathrm{L}_q$}
\put(50,46){\small Reconstruction}
\put(58,42){\small $\mathrm{L}^\epsilon_q$}
\end{picture}
\caption{\label{fig:DOTrecon}Reconstructions of the diffusion coefficient using the suggested equation \eqref{sigma-recon}. $\omega=100 \MHz$ on the left, $\omega=0$ on the right. Original on the top row, non-noisy reconstructions in the middle, noisy reconstrucions on the bottom row. The relative $L^2$ -errors were 19-20\% for all reconstructions.}
\end{figure}

The results of the reconstructions using the truncation above are in figure \ref{fig:DOTrecon}: on top row we have the original diffusion, then non-noisy reconstructions in the middle row, then noisy reconstructions on the bottom row, all using the suggested reconstruction equation \eqref{sigma-recon}. In left column we measure $\Lambda_q$ with $\omega=100 \MHz$ and in right column we use $\omega=0$. The relative $L^2$ -errors were 19-20\% for all reconstructions.

\section{Conclusions}
We developed a numerical method to compute the Faddeev Green's function $g_\lambda(z)$ for negative energy and thus extended the numerical D-bar method from the cases $E=0$ and $E>0$ to negative energy. The computation of the scattering \eqref{scat_matrix} works as expected, the use of noisy DN-matrix results to earlier breakdown as $\abs{\lambda}$ is increased, see figure \ref{scats}. In figure \ref{recon} we see that the D-bar method at negative energy works, but we would hope for a better reconstruction: especially using the scattering from the direct problem, the black line of figure \ref{scats}, we would expect a better reconstruction. This might be improved in the future by improving the computation of $q_0$ from the CGO solution, for example we could use 5-point stencil in \eqref{q0recon-discrete}.

Consequently we found new information on exceptional points for radial potentials, see figures \ref{q1-plane} and \ref{q2-plane}, even for large positive potentials with $\alpha>0$, there are no exceptional points according to our test.

A new method for reconstructing the conductivity at negative energy was presented and tested and found to be working well, see the numerical example \ref{fig:s-recon}. The method works relatively better than the reconstruction of $q_0$ using \eqref{q0recon-discrete} since we don't need the numerical differentiation and the use of ''large $\lambda$''. Despite of additional approximations in assuming that the resulting potential is real and of ''conductivity-type'', our suggested method could be used in DOT to reconstruct the diffusion coefficient, as evidenced by our final result of figure \ref{fig:DOTrecon}.

\section*{Acknowledgements}
This work has been supported by the Academy of Finland (projects 136220, 272803, and 250215 Finnish Centre of Excellence in Inverse Problems Research). JPT was supported in part by the Finnish Cultural Foundation and European Research Council (ERC).

\bibliographystyle{plain}
\bibliography{\refpath Inverse_problems_references}

\end{document}